\documentclass{article}      
\usepackage{amssymb,amsthm}
\title{A Hypothesis on Upper Bound of Goldbach Counting Function}  
\author{Willie B. Wu}      
\date{\today}      

\newtheorem{theorem}{Theorem}[section]
\newtheorem{lemma}[theorem]{Lemma}

\newtheorem{corollary}[theorem]{Corollary}

\newenvironment{definition}[1][Definition]{\begin{trivlist}
\item[\hskip \labelsep {\bfseries #1}]}{\end{trivlist}}
\newenvironment{example}[1][Example]{\begin{trivlist}
\item[\hskip \labelsep {\bfseries #1}]}{\end{trivlist}}

\newenvironment{notations}[1][Notations]{\begin{trivlist}
\item[\hskip \labelsep {\bfseries #1}]}{\end{trivlist}}

\usepackage[top=1in, bottom=1in, left=1in, right=1in]{geometry}
\begin{document}             


\setlength{\parindent}{0pt}

\maketitle                   


\newcounter{N0}
\setcounter{N0}{312}         

\begin{abstract} 
Let $N\geq 1$ be integer, $P\geq 1$ be square-free integer and $S_P(N,x)$ be the number of $n$ between 1 and $x$ such that $(N-n)(N+n)$ is co-prime to $P$. In this paper we propose one hypothesis on upper bound of $S_P(N,x)$ and prove that, under hypothesis $(\ref{upper_bound_hypo_first})$, $S_P(N,N-2)\geq 1$ if $N \geq \arabic{N0}$ and $P$ is the product of all primes $\leq\sqrt{2N}$. Consequently, Goldbach conjecture is true under hypothesis $(\ref{upper_bound_hypo_first})$: there exists $n$ between 1 and $N-2$ such that both $N-n$ and $N+n$ are prime and even number $2N = (N-n)+(N+n)$, sum of two distinct primes if $N \geq \arabic{N0}$. Also, we propose a similar hypothesis on upper bound of $S_P(N,x)$ and prove that, under hypothesis $(\ref{upper_bound_hypo_third})$, the generalized twin prime conjecture is true: for each $N\geq 1$, there are infinitely many pairs of primes $p$ and $q$ such that $q - p = 2N$. 
\end{abstract}

\begin{notations}
We use $\perp$ to indicate two integers are co-prime: $a\perp b$ means $\gcd(a,b) = 1$. $[a,b]$ is the least common multiplier of $a$ and $b$. $p$ and $q$ are primes, $x$ and $y$ are real. $\left\lfloor x\right\rfloor$ is the integer part of $x$ and $\left\{x \right\}$ is the fractional part of $x$. We denote $\mathcal{I}_{x,y}$ to be the set of integers $n$: $x\leq n\leq y$. For any finite set $W$, we denote $|W|$ to be the size of set $W$. We use $\triangleq$ to define new symbols in equation of either side. 
\end{notations} 

\section{Introduction}

Assume $N\geq 1$ is integer and $P\geq 1$ is square-free integer. Goldbach counting function, $S_P(N,x)$ for $x > 0$, is defined to be the number of $n\in \mathcal{I}_{1,x}$ (integers between 1 and $x$) such that $(N-n)(N+n)\perp P$. For application to Goldbach conjecture, we take $P$ to be the product of all primes $\leq \sqrt{2N}$. If $S_P(N,N-2) > 0$, then there is $n\in \mathcal{I}_{1,N-2}$ such that $(N-n)(N+n)\perp P$; it means both $N-n$ and $N+n$ are prime and even number $2N = (N-n)+(N+n)$, sum of two distinct primes. For application to the generalized twin prime conjecture, we take a large $M$ and $P$ to be the product of all primes $\leq M$. If $S_P(N,M^2-N) > 0$, then there is $n\in \mathcal{I}_{1,M^2-N}$ such that $(n-N)(n+N)\perp P$; it means both $q = n+N$ and $p = n-N$ are prime and difference of them, $2N = q-p$, is constant regardless of choice of $M$. If there are infinitely many of such $M$, then there are infinitely many pairs of primes $p$ and $q$ such that $q - p = 2N$. \\ \\
Two basic properties of $S_P(N,x)$ are discussed here: (i) Decomposition of Goldbach counting function: $S_P(N,x)$ is sum of $S_P^d(N,x)$ for all factors $d\mid P$ with respect to $d\perp 6N$, and $S_P(N,x) > 0$ if and only if $S_P^d(N,x) > 0$ for some $d$. (ii) Deduction formula: $S_P(N,x)$ can be expressed in terms of $S_{P'}(N,x)$ for $P' = P/p$ where prime $p\mid P$ and $p\nmid 2N$. Deduction formula for $S_P^d(N,x)$ does exist and a hypothesis is proposed according to the deduction formula for $S_P^p(N,x)$ for prime $p\mid P$ and $p\nmid 6N$. \\ \\
Denote $N_P = \gcd(N,P)$ and $P_n = P/\gcd(P,n)$ for $n\geq 1$. For examples, $P_{6N} = P/\gcd(P,6N)$ and $P_p=P/p$ where $p\mid P$. Let $W_P(N)$ be the set of $n\in \mathcal{I}_{1,P}$ such that $(N-n)(N+n)\perp P$. We will prove 
\begin{eqnarray}
|W_P(N)| = \prod_{p\mid N_P}(p-1) \prod_{p\mid P_{6N}} (p-2) 
\end{eqnarray} 

For $x > 0$, let $S_P(N,x)$ be the number of $n\in \mathcal{I}_{1,x}$ such that $(N-n)(N+n)\perp P$. For $d\mid P_{6N}$, a factor of $P_{6N}$, let $S_P^d(N,x)$ be the number of $n\in \mathcal{I}_{1,x}$ such that $(N-n)(N+n)\perp P$ and $\gcd(P_{6N},n) = d$. Both $S_P(N,x)$ and $S_P^d(N,x)$ are called \emph{Goldbach counting functions}. For $d\mid P_{6N}$, let $W_P^d(N)$ be the set of $n\in W_P(N)$ such that $\gcd(P_{6N},n) = d$. It is clear that $W_P(N)$ is the disjoint union of $W_P^d(N)$ for all $d\mid P_{6N}$. We will prove that size of $W_P^d(N)$ is
\begin{eqnarray}
|W_P^d(N)| = \prod_{p\mid N_P}(p-1) \prod_{p\mid P_{6dN}} (p-3)
\end{eqnarray} 

{\bf UBH:} First upper bound hypothesis on $S_P^p(N,x)$. Assume $N \geq \arabic{N0}$ and $P$ is the product of all primes $\leq \sqrt{2N}$. For any $p\mid P_{6N}$ and any reals $y > x$ and $N/2\leq x < N-1$, the following inequality holds: 
\begin{eqnarray} \label{upper_bound_hypo_first} 
S_P^p(N,y+x) - S_P^p(N,y-x) \leq \frac{3x}{P}|W_P^p(N)|
\end{eqnarray}

This hypothesis says that the number of $n\in \mathcal{I}_{y-x,y+x}$ such that $(N-n)(N+n)\perp P$ and $\gcd(P_{6N},n) = p$ is no more than 150\% of its average. {\bf UBH} $(\ref{upper_bound_hypo_first})$ fails for some small $N$; however, numerical calculation strongly supports this hypothesis for $N\geq 100,000,000$. 
\begin{theorem} \label{Main_Theorem}
Assume $N\geq \arabic{N0}$ and $P$ is the product of all primes $\leq \sqrt{2N}$. If {\bf UBH} $(\ref{upper_bound_hypo_first})$ is true, then 
\begin{eqnarray} 
S_P(N,N-2) \geq 1
\end{eqnarray}
\end{theorem}

One of major steps in the proof of this theorem is deduction formula $(\ref{deduction_form})$ for Goldbach counting function. It is a formula for $S_P(N,x)$ in terms of $S_{P_p}(N,x)$ for $p\mid P_{2N}$. For prime $p\mid P$, let $\overline{P_p}$ be an inverse of $P_p$ in $p$, satisfying $\overline{P_p} P_p \equiv 1$ mod $p$, and $\bar p$ be an inverse of $p$ in $P_p$, satisfying $\bar pp \equiv 1$ mod $P_p$. 

\begin{theorem} {\bf (Deduction formula for $S_P(N, x)$)} 
For $p\mid P_{2N}$, if $x$ is not integer and $x < N\overline{P_p}P_p$, then 
\begin{eqnarray} \label{deduction_form}
S_P(N,x) = S_{P_p}(N,x) - S_{P_p}\left(\bar pN,\frac{N\overline{P_p}P_p+x}{p}\right) + S_{P_p}\left(\bar pN,\frac{N\overline{P_p}P_p-x}{p}\right)
\end{eqnarray}
\end{theorem}
\begin{proof}
Let $m = N\overline{P_p}P_p$. We only need to prove for $x = n+\frac{1}{2}$ where $n = 0, 1, \cdots, m-1$. For $n=0$, we have $S_{P_p}(N,x) = S_P(N,x) = 0$ and $S_{P_p}\left(\bar pN,\frac{m+x}{p}\right) = S_{P_p}\left(\bar pN,\frac{m-x}{p}\right)$. For $0<n<m$, we need to prove 
\begin{eqnarray}
S_{P_p}(N,x) - S_P(N,x) = S_{P_p}\left(\bar pN,\frac{m+x}{p}\right) - S_{P_p}\left(\bar pN,\frac{m-x}{p}\right)
\end{eqnarray}
as $x$ changes from $x = n-\frac{1}{2}$ to $x = n+\frac{1}{2}$. Let us look at one case: $S_{P_p}\left(\bar pN,\frac{m+x}{p}\right)$ increases by 1 from $x = n-\frac{1}{2}$ to $x = n+\frac{1}{2}$. In this case, $n' = \frac{m+x'}{p}$ is integer for some $x'$ and $(\bar pN+n')(\bar pN-n')\perp P_p$. Since $\frac{m+x'}{p}$ is integer, then $x' = n$. $p\mid m + x' = N\overline{P_p}P_p + n$ means $p\mid N + n$. Since $(\bar pN+n')(\bar pN-n')\perp P_p$, then $(N+pn')(N-pn')\perp P_p$ and $(N+n)(N-n)\perp P_p$. Thus, $S_{P_p}(N,x) - S_P(N,x)$ increases by 1 from $x = n-\frac{1}{2}$ to $x = n+\frac{1}{2}$ and $S_{P_p}\left(\bar pN,\frac{m-x}{p}\right)$ has no change. The rest of proof is to verify the other cases. 
\end{proof}

It is not hard to check the other cases; however, we skip the verification here since we will give an ``analytic" proof after the Goldbach cosine sum-product formula is established. For $p\mid P$, we define 
\begin{eqnarray}
\alpha_p(P,m) = \left\{\begin{array}{ll}
p-2 & \mbox{if } p\mid m, \\
-2\cos\frac{2m\overline{P_p}\pi}{p} & \mbox{if } p\nmid m
\end{array}\right.
\end{eqnarray}
We are going to show the \emph{Goldbach cosine sum-product formula} over $W_P(N)$ for integer $k$:
\begin{eqnarray} 
C_P(N,k) \triangleq \sum_{n\in W_P(N)} \cos\frac{2nk\pi}{P} = \mu(N_P)\prod_{p\mid\gcd(k,N_P)}(1-p) \prod_{p\mid P_{6N}} \alpha_p(P,kN)
\end{eqnarray}
where $\mu$ is the M\"{o}bius function, and the \emph{Goldbach cosine sum-product formula} over $W_P^d(N)$ for $d\mid P_{6N}$: 
\begin{eqnarray}
C_P^d(N,k) \triangleq \sum_{n\in W_P^d(N)} \cos\frac{2nk\pi}{P} = \mu(N_P)\prod_{p\mid\gcd(k,N_P)}(1-p) \prod_{p\mid P_{6dN}} (\alpha_p(P,kN)-1)
\end{eqnarray}

Next we will prove the Goldbach counting formula when $x$ is not integer: 
\begin{eqnarray}
S_P(N,x) = |W_P(N)|\frac{x}{P} - t_P(N) - T_P(N,x) 
\end{eqnarray}
where
\begin{eqnarray}
t_P(N) = \left\{\begin{array}{cl} 
\frac{1}{2} & \mbox{if } N_P=1, \\ 0 & \mbox{if } N_P>1
\end{array}\right. \mbox{ and } T_P(N,x) = - \sum_{k = 1}^\infty \frac{C_P(N,k)}{k\pi} \sin \frac{2k\pi x}{P}
\end{eqnarray}

By Goldbach cosine sum-product formula over $W_P(N)$, we can get the deduction formula for $C_P(N,k)$: 
\begin{eqnarray}
C_P(N,k) = C_{P_p}(\bar pN,k) \alpha_p(P,kN)
\end{eqnarray}
for $p\mid P_{2N}$. By this formula, we are able to show    
\begin{theorem} {\bf (Deduction formula for $T_P(N, x)$)} 
For $p\mid P_{2N}$, if $x$ is not integer and $x < N\overline{P_p}P_p$, then 
\begin{eqnarray} 
T_P(N,x) = T_{P_p}(N,x) - T_{P_p}\left(\bar pN,\frac{N\overline{P_p}P_p+x}{p}\right) + T_{P_p}\left(\bar pN,\frac{N\overline{P_p}P_p-x}{p}\right)
\end{eqnarray}
\end{theorem}

Deduction formula for $S_P(N, x)$ can be derived from deduction formula for $T_P(N, x)$. $T_P(N,x)$ is called \emph{error} term; however, it is not small in general. We will see $T_P(N,N)/S_P(N,N) \approx 0.26$ by numerical calculations for large $N$. For $d\mid P_{6N}$, we define another error term: 
\begin{eqnarray}
T_P^d(N,x) = -\sum_{k = 1}^\infty \frac{C_P^d(N,k)}{k\pi} \sin \frac{2k\pi x}{P}
\end{eqnarray}
and will prove that
\begin{eqnarray}
S_P^d(N,x) = |W_P^d(N)|\frac{x}{P} - t_P^d(N) - T_P^d(N,x) 
\end{eqnarray}
where $t_P^d(N) = 0$ if $d<P_{6N}$ and $t_P^d(N) = t_P(N)$ if $d=P_{6N}$. By this formula, {\bf UBH} $(\ref{upper_bound_hypo_first})$ can be given equivalently as, for $p\mid P_{6N}$, 
\begin{eqnarray} \label{upper_bound_hypo_second} 
T_P^p(N,y-x) - T_P^p(N,y+x) \leq \frac{x}{P}|W_P^p(N)| 
\end{eqnarray}

We will show $C_P^d(N,k) = C_P^{dp}(N,k)(\alpha_p(P,kN)-1)$ for $d\mid P_{6N}$ and $p\mid P_{6dN}$, and  
\begin{eqnarray}
T_P^d(N,x) = T_{P_p}^d(N,x) - T_P^{dp}(N,x) - T_P^{dp}(N,N\overline{P_p}P_p+x) + T_P^{dp}(N,N\overline{P_p}P_p-x) 
\end{eqnarray}

By taking $d = 1$, we get the following after {\bf UBH} $(\ref{upper_bound_hypo_second})$ is applied with $y = N\overline{P_p}P_p$: 
\begin{eqnarray}
T_P^1(N,x) \leq T_{P_p}^1(N,x) - T_P^p(N,x) + \frac{x}{P}|W_P^p(N)| 
\end{eqnarray}

for $N/2\leq x<N-1$. If $T_P^p(N,x) < \frac{x}{P}|W_P^p(N)|$, then $S_P^p(N,x) > 0$ and $S_P(N,N-2)\geq 1$. Otherwise, we have $T_P^p(N,x) \geq \frac{x}{P}|W_P^p(N)|$ and $T_P^1(N,x) \leq T_{P_p}^1(N,x)$. It is the time to present the following: 

\begin{theorem} 
Assume $x > 0$. If $T_P^1(N,x) \leq T_{P_p}^1(N,x)$ for each $p\mid P_{6N}$, then $T_P^1(N,x) \leq T_{P_d}^1(N,x)+t_{P_d}^1(N)$ for any $d\mid P_{6N}$. 
\end{theorem}
It turns out the proof of this theorem is quite simple. Now we have a small error term $T_P^1(N,x) \leq T_{P_d}^1(N,x)+t_P^d(N)$ by selecting $d = P_{6N}$ and $N/2 \leq x < N-1$. This is the outline of our major steps to prove Theorem \ref{Main_Theorem}. Similarly, we will prove the generalized twin prime conjecture under {\bf UBH$'$}: \\ \\ 
{\bf UBH$'$:} Second upper bound hypothesis on $S_P^p(N,x)$. For given $N\geq 1$, there are infinitely many integers $M \geq 2N+1$ such that for $M^2-7N \leq x < M^2-N < y$ and for each $p\mid P_{6N}$ where $P$ is the product of all primes $\leq M$, the following inequality holds: 
\begin{eqnarray} \label{upper_bound_hypo_third} 
S_P^p(N,y+x) - S_P^p(N,y-x) \leq \frac{3x}{P}|W_P^p(N)|
\end{eqnarray}

\begin{theorem} 
If {\bf UBH$'$} $(\ref{upper_bound_hypo_third})$ is true for given $N\geq 1$, then there are infinitely many pairs of primes $p$ and $q$ such that $q - p = 2N$. 
\end{theorem}

Let's start with the decomposition of set $W_P(N)$.


\section{Decomposition of $W_P(N)$}

There are several ways to obtain the formula for $|W_P(N)|$. By use of Chinese remainder theorem, we get it easily. Here it is. 
\begin{theorem}
Size of $W_P(N)$ is 
\begin{eqnarray}
|W_P(N)| = \prod_{p\mid N_P}(p-1) \prod_{p\mid P_{2N}} (p-2) = \prod_{p\mid N_P}(p-1) \prod_{p\mid P_{6N}} (p-2) 
\end{eqnarray} 
\end{theorem}
\begin{proof}
If $2\mid P_N$, let $\mathcal{J}_2 = \{0\}$. For $p\mid P_{2N}$, let $\mathcal{J}_p = \mathcal{I}_{0,p-1}\setminus\{a,b\}$ where $a$ and $b\in \mathcal{I}_{1,p-1}$ are the solutions $s$ of $p\mid s + N$ and $p\mid s - N$ respectively. Since $a\ne b$, then $|\mathcal{J}_p|=p-2$. For $n\in W_P(N)$ and $p\mid P$, let $h_p\equiv n$ mod $p$ and $h_p\in \mathcal{I}_{0,p-1}$. We are going to prove $h_p\in \mathcal{I}_{1,p-1}$ if $p\mid N_P$ and $h_p\in \mathcal{J}_p$ if $p\mid P_N$. First we assume $p\mid N_P$, then $p\nmid (N-n)(N+n)$ means $p\nmid n$ and $h_p\in \mathcal{I}_{1,p-1}$. Next we assume $p\mid P_N$. If $p = 2$, then $2\nmid N$ and $2\nmid(N-n)(N+n)$ mean $n$ is even and $h_2 = 0$. If $p>2$, then $p\nmid (N-n)(N+n)$ means $h_p\not\equiv \pm N$ mod $p$ and therefore, $h_p\in \mathcal{J}_p$. If $n'\in W_P(N)$ and $p\mid n-n'$ for all $p\mid P$, then $n'=n$ and 
\begin{eqnarray}
|W_P(N)| \leq \prod_{p\mid N_P}|\mathcal{I}_{1,p-1}| \prod_{p\mid P_N} |\mathcal{J}_p| = \prod_{p\mid N_P}(p-1) \prod_{p\mid P_{2N}} (p-2) 
\end{eqnarray}
Here we have $|\mathcal{J}_2|=1$ if $2\mid P_N$. Now we pick one value $f_p\in \mathcal{I}_{1,p-1}$ for $p\mid N_P$ and one value $f_p\in \mathcal{J}_p$ for $p\mid P_N$, then the system of equations 
\begin{eqnarray}
z \equiv f_p \mbox{ mod } p \mbox{ for all } p\mid P
\end{eqnarray}
has one solution $n$ between 1 and $P$ by Chinese remainder theorem. Thus, $z = n\in W_P(N)$ and 
\begin{eqnarray}
|W_P(N)| \geq \prod_{p\mid N_P}(p-1) \prod_{p\mid P_{2N}} (p-2) = \prod_{p\mid N_P}(p-1) \prod_{p\mid P_{6N}} (p-2) 
\end{eqnarray} 
Here we have $|\mathcal{J}_3| = 1$ if $3\mid P_{2N}$. That completes the proof.
\end{proof}

\begin{definition}
For $d\mid P_{6N}$, a factor of $P_{6N}$, let $W_P^d(N)$ be the set of $n\in W_P(N)$ such that $\gcd(P_{6N},n) = d$. 
\end{definition} 

\begin{theorem}
{\bf (Decomposition of $W_P(N)$)} $W_P(N)$ is the disjoint union of $W_P^d(N)$ for all $d\mid P_{6N}$.
\end{theorem}
\begin{proof}
For $n\in W_P(N)$, let $d=\gcd(P_{6N}, n)$, then $n\in W_P^d(N)$ and $W_P(N)$ is the union of $W_P^d(N)$ for all $d\mid P_{6N}$. If $n\in W_P^d(N)\cap W_P^{d'}(N)$, then $d = \gcd(P_{6N}, n) = d'$ and $W_P^d(N)$ is distinct for each $d\mid P_{6N}$. 
\end{proof}

We define $c = \frac{\gcd(NP,6)}{\gcd(N,6)}$, the \emph{index} of $(P,N)$. Clearly $c\mid 6$. 
\begin{lemma}
{\rm i.} $c\mid n$ if $n\in W_P(N)$. {\rm ii.} $c\perp N$. {\rm iii.} $P_{cN} = P_{6N}$. 
\end{lemma}
\begin{proof}
Assume $n\in W_P(N)$, then $P\perp (N-n)(N+n)$. If $2 \mid c$, then, by definition of $c$, $N$ is odd, $P$ is even and $2\mid n$ since $2\nmid (N-n)$. If $3\mid c$, then $3\mid P$, $3\nmid N$ and $3\mid n$ since $3\nmid (N-n)(N+n)$. Thus, part i and part ii are proved. For part iii, we only need to prove $2\nmid P_{cN}$ and $3\nmid P_{cN}$; it is obvious by definition of $c$. 
\end{proof} 

\begin{definition}
For $d\mid P_{6N}$, let $P_{cd}^\perp$ be the set of $k\perp P_{cd}$ and $1\leq k\leq P_{cd}$. 
\end{definition} 

\begin{theorem}
Assume $n\in W_P$ and $d\mid P_{6N}$. $n\in W_P^d$ if and only if $n=cdk$ and $k\in P_{cd}^\perp$. 
\end{theorem}
\begin{proof}
Assume $n\in W_P^d$, then $d = \gcd(P_{6N},n)$. Since $cd \mid n \leq P$ then $n=cdk$ for some $k\leq P_{cd}$. 
Since $n\in W_P$, then $n\perp N_P$ and $k\perp N_P$. Since $\gcd(P_{cN}, cdk) = d$, then $\gcd(P_{cdN}, k) = 1$ and $k\perp P_{cdN}$. Thus, $k\perp N_PP_{cdN}=P_{cd}$ and $k\in P_{cd}^\perp$. Now assume $n=cdk$ and $k\in P_{cd}^\perp$. Let $d'= \gcd(P_{6N},n) = \gcd(P_{cN},cdk)$, then $d'/d = \gcd(P_{cdN},ck) = \gcd(P_{cdN},k)$. Since $\gcd(P_{cd},k) = 1$, then $\gcd(P_{cdN},k) = 1$ and $d'=d$. 
\end{proof}

For $d\mid P_{6N}$, let $V_P^d(N)$ be the set of $k\in P_{cd}^\perp$ such that $cdk\in W_P^d(N)$. Clearly $|V_P^d(N)| = |W_P^d(N)|$. 

\begin{theorem}
For $d\mid P_{6N}$, size of $V_P^d(N)$ is
\begin{eqnarray}
|V_P^d(N)| = \prod_{p\mid N_P}(p-1) \prod_{p\mid P_{6dN}} (p-3) 
\end{eqnarray}
\end{theorem}
\begin{proof}
For $p\mid P_{6dN}$, let $\mathcal{K}_p = \mathcal{I}_{1,p-1} \setminus \{a,b\}$ where $a$ and $b\in \mathcal{I}_{1,p-1}$ are the solutions $s$ of $p\mid cds+N$ and $p\mid cds-N$ respectively. Since $a\ne b$, then $|\mathcal{K}_p| = p-3$. For $k\in V_P^d(N)$, let $h_p\equiv k$ mod $p$ and $h_p\in \mathcal{I}_{0,p-1}$ for $p\mid P$. Since $V_P^d(N)\subset P_{cd}^\perp$, then $h_p\in \mathcal{I}_{1,p-1}$. For $p\mid P_{6dN}$, we have $p\nmid (N-cdk)(N+cdk)$; thus, $cdh_p\not\equiv \pm N$ mod $p$ and $h_p\in \mathcal{K}_p$. Therefore, 
\begin{eqnarray}
|V_P^d(N)| \leq \prod_{p\mid N_P}|\mathcal{I}_{1,p-1}| \prod_{p\mid P_{6dN}} |\mathcal{K}_p| = \prod_{p\mid N_P}(p-1) \prod_{p\mid P_{6dN}} (p-3) 
\end{eqnarray}
Now we pick one value $f_p\in \mathcal{I}_{1,p-1}$ for $p\mid N_P$ and one value $f_p\in \mathcal{K}_p$ for $p\mid P_{6dN}$. Since $P_{cd} = N_PP_{6dN}$, then the system of equations 
\begin{eqnarray}
z \equiv f_p \mbox{ mod } p \mbox{ for all } p \mid P_{cd}
\end{eqnarray}
has one solution $k$ between 1 and $P_{cd}$ by Chinese remainder theorem. Thus, $z = k\in V_P^d(N)$ and the size of $V_P^d(N)$ is
\begin{eqnarray}
|V_P^d(N)| \geq \prod_{p\mid N_P}|\mathcal{I}_{1,p-1}| \prod_{p\mid P_{6dN}} |\mathcal{K}_p| = \prod_{p\mid N_P}(p-1) \prod_{p\mid P_{6dN}} (p-3) 
\end{eqnarray}
That completes the proof. 
\end{proof}

\begin{theorem}
Let $K$ be a square-free integer and let $h_p$ be given for prime $p\mid K$, then 
\begin{eqnarray}
\prod_{p\mid K} (h_p+1) = \sum_{d\mid K} \prod_{p\mid K_d} h_p 
\end{eqnarray}
where $d$ goes over all factors of $K$ and $K_d = K/d$. 
\end{theorem}

An easy way to understand this formula is to treat each $h_p$ as an indeterminate in polynomial. This formula will be applied to the proof of Goldbach cosine sum-product formula over $W_P^d(N)$ and two more places: one is in the following example and one in the proof of Goldbach momentum formula.  

\begin{example}
We take $K = P_{6N}$ and $h_p = p-3$ for $p\mid P_{6N}$, then 
\begin{eqnarray}
\prod_{p\mid P_{6N}} (p-2) & = & \sum_{d\mid P_{6N}} \prod_{p\mid P_{6dN}} (p-3) \\
\prod_{p\mid N_P}(p-1)\prod_{p\mid P_{6N}} (p-2) & = & \sum_{d\mid P_{6N}} \prod_{p\mid N_P}(p-1) \prod_{p\mid P_{6dN}} (p-3) \\
|W_P(N)| & = & \sum_{d\mid P_{6N}} |W_P^d(N)| 
\end{eqnarray}
\end{example}

\begin{theorem}
{\bf (Goldbach momentum formula)} For any $f_p$ defined for each $p\mid P$, 
\begin{eqnarray} 
\sum_{n\in W_P(N)} \prod_{p\mid P_n} \frac{f_p-1}{p-1} = \prod_{p\mid N_P}(f_p-1) \prod_{p\mid P_{6N}} \left(f_p - 2\frac{f_p-1}{p-1}\right) 
\end{eqnarray}
\end{theorem}
\begin{proof}
For $p\mid P$ and given $f_p$, we define 
\begin{eqnarray}
h_p = \frac{f_p-1}{p-1}(p-3) = f_p - 2\frac{f_p-1}{p-1} - 1
\end{eqnarray}
Thus, 
\begin{eqnarray}
\prod_{p\mid P_{6N}}\left(f_p - 2\frac{f_p-1}{p-1}\right) = \prod_{p\mid P_{6N}} (h_p+1)  & = & \sum_{d\mid P_{6N}} \prod_{p\mid P_{6dN}} h_p = \sum_{d\mid P_{6N}} \prod_{p\mid P_{6dN}} \frac{f_p-1}{p-1}(p-3) \\
\prod_{p\mid N_P}(p-1) \prod_{p\mid P_{6N}}\left(f_p - 2\frac{f_p-1}{p-1}\right) & = & \sum_{d\mid P_{6N}} \prod_{p\mid N_P}(p-1) \prod_{p\mid P_{6dN}} \frac{f_p-1}{p-1}(p-3) \\
& = & \sum_{d\mid P_{6N}} |V_P^d(N)| \prod_{p\mid P_{6dN}} \frac{f_p-1}{p-1}
\end{eqnarray}
Since $P_{6dN} = P_{cdkN}$ for any $d\mid P_{6N}$ and $k\in V_P^d(N)$, then 
\begin{eqnarray}
\prod_{p\mid N_P}(p-1) \prod_{p\mid P_{6N}}\left(f_p - 2\frac{f_p-1}{p-1}\right) = \sum_{d\mid P_{6N}} \sum_{k\in V_P^d(N)} \prod_{p\mid P_{cdkN}} \frac{f_p-1}{p-1} 
= \sum_{n\in W_P(N)} \prod_{p\mid P_{nN}} \frac{f_p-1}{p-1} 
\end{eqnarray}
For $n\in W_P(N)$, we have $n\perp N_P$ and $P_n = N_PP_{nN}$; thus, 
\begin{eqnarray}
\prod_{p\mid P_n} \frac{f_p-1}{p-1} & = & \prod_{p\mid N_P} \frac{f_p-1}{p-1}\prod_{p\mid P_{nN}} \frac{f_p-1}{p-1} \\
\sum_{n\in W_P(N)} \prod_{p\mid P_n} \frac{f_p-1}{p-1} & = & \prod_{p\mid N_P} \frac{f_p-1}{p-1} \sum_{n\in W_P(N)} \prod_{p\mid P_{nN}} \frac{f_p-1}{p-1} \\
& = & \prod_{p\mid N_P} \frac{f_p-1}{p-1} \prod_{p\mid N_P}(p-1) \prod_{p\mid P_{6N}}\left(f_p - 2\frac{f_p-1}{p-1}\right)
\end{eqnarray}
and Goldbach momentum formula follows. 
\end{proof}

Goldbach momentum formula was first discovered for $f_p = p^s$ by another method \cite{Wu}. We are going to prove the Goldbach cosine sum-product formula in the next several sections.

\section{Modulo Set over Square-Free Integer} 

Let us start with the following theorem. 
\begin{theorem}
Let $N\geq 1$ be an integer and $P\geq 1$ be a square-free integer. For $a\mid P$ and $b\mid P$, there exists an integer $n$ such that $a\mid N-n$ and $b\mid N+n$ if and only if $\gcd(a, b)\mid 2N$. 
\end{theorem} 
\begin{proof}
If there is $n$ such that $a\mid N-n$ and $b\mid N+n$, then $n = N - ak$ for some $k$ and $n = bl - N$ for some $l$. Thus, $n = N - ak = bl - N$ and $2N = ak+bl$. Therefore, $\gcd(a,b)\mid 2N$. Now assume $\gcd(a,b)\mid 2N$. Let $n_1$ be the least integer $n\geq 0$ such that $b\mid N+n$ and $N_1 = (N + n_1)/b$. If $b\mid N+n$, then $n = n_1 + bk$ for some $k$ and 
\begin{eqnarray}
N-n = N- n_1 - bk = 2N - (N+n_1) - bk = 2N - bN_1 - bk
\end{eqnarray}
Let $d = \gcd(a, b)$, $a' = a / d$ and $b' = b/d$, then $a'\perp b'$. Since $d\mid 2N$, then $a\mid N-n$ is equivalent to $a'\mid N_2-b'N_1-b'k$ where $N_2 = 2N/d$. Let $k_1$ be the least integer $k \geq 0$ such that $a'\mid N_2-b'N_1 -b'k$ (the existence of such $k$ is due to $a'\perp b'$), then $n = n_1 + bk_1$ meets the requirement. 
\end{proof}

\begin{definition}
Let $N\geq 1$ be an integer and $P\geq 1$ be a square-free integer. Let $Q$ be the set of pairs $(a,b)$ of positive integers such that $[a,b]\mid P$ and $\gcd(a,b)\mid 2N$. For $(a, b)\in Q$, let $m_{ab}$ be the least integer $n\geq 0$ such that $a\mid N-n$ and $b\mid N+n$. We call $\{m_{ab}\}$ the \emph{modulo set} \index{modulo set} in respect of $N$ and $P$. 
\end{definition}

We assume $N$ and $P$ are given throughout this paper; therefore so is $Q$. The following three theorems uncover some properties on the modulo set $\{m_{ab}\}$. 
\begin{theorem}
If $(a,b),(a',b)\in Q$ and $a\mid a'$, then $m_{a'b} = m_{ab}+k[a,b]$ for some $k\geq 0$. Likewise, if $(a,b),(a,b')\in Q$ and $b\mid b'$, then $m_{ab'} = m_{ab}+k'[a,b]$ for some $k'\geq 0$. 
\end{theorem}
\begin{proof}
Since $a\mid a' \mid N-m_{a'b}$ and $a\mid N-m_{ab}$, then $a\mid m_{a'b} - m_{ab}$. Similarly, $b\mid m_{a'b} - m_{ab}$ and $[a,b]\mid m_{a'b} - m_{ab}$. Thus, $m_{a'b} = m_{ab}+k[a,b]$ and $k\geq 0$ since $m_{a'b} \geq m_{ab}$. 
\end{proof}

\begin{theorem}
Assume $(a,b)\in Q$ and prime $p\mid \gcd(2N,b)$. Let 
\begin{eqnarray}
a' = \left\{\begin{array}{ll}
a/p & \mbox{\rm if } p\mid a, \\ ap & \mbox{\rm if } p\nmid a
\end{array}\right.
\end{eqnarray}
then $(a',b)\in Q$ and $m_{a'b} = m_{ab}$. 
\end{theorem}
\begin{proof}
First, $\gcd(a',b) = \gcd(a,b)$ and $(a',b)\in Q$. Next we assume $p = 2$. Since $2\mid b\mid N+m_{ab}$, then $2\mid N-m_{ab}$ and $a'\mid N-m_{ab}$. Now we assume $p$ is odd and $p\mid N$. Since $p\mid b\mid N+m_{ab}$, then $p\mid m_{ab}$, $p\mid N-m_{ab}$ and $a'\mid N-m_{ab}$. In either case we have $m_{a'b} \leq m_{ab}$. If $a' = ap$ then we have $m_{a'b} = m_{ab}$. If $a' = a/p$, then $a'' \triangleq (a')'$ and $m_{a''b} = m_{a'b}$. Since $a = a''$, then the proof is complete. 
\end{proof}

\begin{theorem}
For $(a,b)\in Q$, we have
\begin{eqnarray}
m_{ab} + m_{ba} = \left\{\begin{array}{ll}
0 & \mbox{\rm if } [a,b] \mid N, \\
{[a,b]} & \mbox{\rm if } [a,b] \nmid N \end{array}\right. 
\end{eqnarray}
\end{theorem}
\begin{proof}
If $[a,b] \mid N$, then it is clear that $m_{ab} = m_{ba} = 0$. If $[a,b] \nmid N$, then $m_{ab} > 0$. Let $n = [a,b] - m_{ab}$, then $n > 0$, $a\mid N + n$ and $b\mid N - n$; thus, $m_{ba} = n$. 
\end{proof}

\begin{definition}
For $n\geq 1$, let $U_n$, the \emph{unit set} at $n$, be the set of $(a,b)\in Q$ such that $[a,b] \mid n - m_{ab}$, let $U_n^\ast$ be the set of $(a,b)\in Q$ such that $[a,b]\mid n + m_{ab}$, called the \emph{dual} of unit set $U_n$. 
\end{definition}

It is clear that $(a,b)\in U_n$ if and only if $(b,a)\in U_n^\ast$. M\"{o}bius function $\mu$ is widely used in the classic sieve method \cite{Charles}. For $n\geq 1$, we define the \emph{unit value} at $n$ as 
\begin{eqnarray}
u_n = \sum_{(a,b)\in U_n} \mu(a)\mu(b)
\end{eqnarray}
\begin{theorem} 
For $n\geq 1$, the unit value $u_n$ is given as follows: 
\begin{eqnarray} 
u_n = \left\{\begin{array}{ll}
1 & \mbox{\rm if } (N-n)(N+n)\perp P, \\ 0 & \mbox{\rm otherwise}
\end{array}\right.
\end{eqnarray}
\end{theorem}
\begin{proof}
First, we define $a_n = \gcd(N-n,P)$ and $b_n = \gcd(N+n,P)$. If $(a,b)\in U_n$, then $n = m_{ab} + k[a,b]$ for some $k$. Since $a\mid N-m_{ab}$ and $a\mid n - m_{ab}$, then $a\mid N-n$, therefore $a\mid a_n$; similarly, $b\mid b_n$. Conversely, if $a\mid a_n$ and $b\mid b_n$, then $a\mid N-n$ and $b\mid N+n$. Thus, $n = m_{ab}+k[a,b]$ for some $k$ and $(a,b)\in U_n$. Thus, $U_n$ is the set of $(a,b)\in Q$ such that $a\mid a_n$ and $b\mid b_n$. Now we have
\begin{eqnarray}
u_n = \sum_{(a,b)\in U_n} \mu(a)\mu(b) = \sum_{a\mid a_n} \mu(a) \sum_{b\mid b_n} \mu(b)
\end{eqnarray}
Since $a_n = b_n = 1$ if and only if $(N-n)(N+n)\perp P$, then the theorem follows. 
\end{proof} 

By use of $u_n$, we have for integer $k$, 
\begin{eqnarray}
\sum_{n\in W_P(N)} \cos\frac{2nk\pi}{P} = \sum_{n = 1}^P u_n\cos\frac{2nk\pi}{P}
\end{eqnarray}

\section{Sum-Sieve Equation}

For $(a,b)\in Q$ and $x\geq 0$, we define 
\begin{eqnarray}
n_{ab}(x) = \left\lfloor\frac{x+m_{ab}}{[a,b]}\right\rfloor = \frac{x+m_{ab}}{[a,b]} - \left\{\frac{x+m_{ab}}{[a,b]}\right\}
\end{eqnarray}
the integer part of $\frac{x+m_{ab}}{[a,b]}$, where $\left\{y \right\}$ is the fractional part of $y$. For integer $n \geq 1$, let 
\begin{eqnarray}
n_{ab}(n^-) = \lim_{x\rightarrow n^-}n_{ab}(x)
\end{eqnarray}

Since $n_{ab}(n) \ne n_{ab}(n^-)$ if and only if $[a,b]\mid n+m_{ab}$, then we have the following lemma: 
\begin{lemma} \label{U_n_ab_lemma}
For $n \geq 1$, $n_{ab}(n) \ne n_{ab}(n^-)$ if and only if $(a,b)\in U_n^\ast$. 
\end{lemma}

For $(a,b)\in Q$ and $x\geq 0$, we define 
\begin{eqnarray}
A_{ab} = [a,b] - 2m_{ab} \mbox{ and } B_{ab}(x) = A_{ab} + 2n_{ab}(x)[a,b] 
\end{eqnarray}
\begin{theorem}
{\bf (Sum-sieve equation)} Assume $N$, $P$ and $Q$ are given. For $x\geq 0$ and real $s$ such that $\sin Ps \ne 0$, we have
\begin{eqnarray} 
E(x, s) \triangleq \sum_{(a,b)\in Q}\mu(a)\mu(b) \frac{e^{iB_{ab}(x)s} - e^{iA_{ab}s}}{2i\sin [a,b]s} = \sum_{1\leq n \leq x} u_ne^{2ins} 
\end{eqnarray}
\end{theorem}
where $i = \sqrt{-1}$ is the imaginary unit of complex number.
\begin{proof}
The equation is true for $x = 0$. $E(x, s)$, as a function of $x$, is a step function jumping probably at $x = n$, the positive integers. For $0 < n \leq x$, if $(a,b)\in U_n^\ast$, then $n_{ab}(n^-) = n_{ab}(n) - 1$, $B_{ab}(n) = 2n + [a,b]$ and $B_{ab}(n^-) = 2n - [a,b]$. By Lemma $\ref{U_n_ab_lemma}$, we have
\begin{eqnarray}
E(n, s) - E(n^-, s) & = & \sum_{(a,b)\in U_n^\ast}\mu(a)\mu(b) \frac{e^{iB_{ab}(n)s} - e^{iB_{ab}(n^-)s}}{2i\sin [a,b]s} \\
& = & \sum_{(a,b)\in U_n^\ast}\mu(a)\mu(b) \frac{e^{i(2n + [a,b])s} - e^{i(2n - [a,b])s}}{2i\sin [a,b]s} \\ 
& = & e^{2ins} \sum_{(a,b)\in U_n^\ast}\mu(a)\mu(b) \frac{e^{i[a,b]s} - e^{-i[a,b]s}}{2i\sin [a,b]s} = u_n e^{2ins} 
\end{eqnarray} 
Thus, step function $E(x, s)$ is the sum of all these terms. 
\end{proof}

By taking the real and imaginary parts on $E(x,s)$ respectively, we get
\begin{eqnarray}
E_1(x,s) & \triangleq & \sum_{(a,b)\in Q}\mu(a)\mu(b) \frac{\sin B_{ab}(x)s - \sin A_{ab}s}{2\sin [a,b]s} = \sum_{1\leq n \leq x} u_n \cos 2ns \\
E_2(x,s) & \triangleq & \sum_{(a,b)\in Q}\mu(a)\mu(b) \frac{\cos B_{ab}(x)s - \cos A_{ab}s}{2\sin [a,b]s} = -\sum_{1\leq n \leq x} u_n \sin 2ns 
\end{eqnarray}

$E_1(x,s)$ and $E_2(x,s)$ are called cosine and sine formula on modulo set respectively. 
\begin{definition}
For any factor $d\mid P$, let $Q_d$ be the set of $(a,b)\in Q$ such that $ab=P_d$ and $b\perp 2N$.
\end{definition}

\begin{lemma}
For any function $F(n,m)$ defined on integers $n$ and $m$, we have
\begin{eqnarray}
\sum_{(a,b)\in Q}\mu(a)\mu(b)F(m_{ab}, [a,b]) = \sum_{d\mid P} \mu(P_d) \sum_{(a,b)\in Q_d}F(m_{ab}, [a,b])
\end{eqnarray}
\end{lemma}
\begin{proof}
For $d\mid P$ and $(a,b)\in Q_d$, $\mu(a)\mu(b) = \mu(P_d)$. Let $Q_d'$ be the set of $(a,b)\in Q\setminus Q_d$ such that $[a,b]=P_d$, then
\begin{eqnarray}
\sum_{(a,b)\in Q}\mu(a)\mu(b)F(m_{ab}, [a,b]) = \sum_{d\mid P} \mu(P_d)\sum_{(a,b)\in Q_d}F(m_{ab}, [a,b]) + \sum_{d\mid P} \sum_{(a,b)\in Q_d'}\mu(a)\mu(b)F(m_{ab}, [a,b])
\end{eqnarray} 
We need to prove the second term is 0. Let $P_d = p_1\cdots p_nq_1\cdots q_m$ be all the prime factors such that $p_1\cdots p_n\mid 2N$ and $q_1\cdots q_m \perp 2N$. Since $(a,b)\in Q_d$ if and only if $b\mid q_1\cdots q_m$ and $a = P_d/b$, then $(a,b)\mid Q_d'$ if and only if $p_j\mid b$ for some $j$ and $[a,b]=P_d$. For $(a,b)\in Q_d'$, let $l$ be the maximal index of $j$ such that $p_j \mid b$. Therefore, $p_l \mid \gcd(2N,b)$; we define $a' = a / p_l$ if $p_l \mid a$ and $a' = ap_l$ if $p_l\nmid a$. Thus, $(a',b)\mid Q_d'$, $\mu(a) = -\mu(a')$ and $m_{ab} = m_{a'b}$. Now we have 
\begin{eqnarray} 
2\sum_{(a,b)\in Q_d'}\mu(a)\mu(b) F(m_{ab},[a,b])= \sum_{(a,b)\in Q_d'}(\mu(a)+\mu(a'))\mu(b) F(m_{ab},[a,b]) = 0 
\end{eqnarray}
This completes the proof. 
\end{proof}

By application of this lemma to $E_1(x,s)$, we have 
\begin{eqnarray}
E_1(x,s) = \sum_{1\leq n \leq x} u_n\cos 2ns = \sum_{d\mid P} \mu(P_d) \sum_{(a,b)\in Q_d} 
\frac{\sin B_{ab}(x)s - \sin A_{ab}s}{2\sin P_ds} 
\end{eqnarray}

\section{Goldbach Cosine Sum-Product Formula} 

First, we present the following cosine formula: 
\begin{theorem}
For integer $k\geq 1$, 
\begin{eqnarray} \label{cosine_k_sum_sieve_eq}
\sum_{n\in W_P(N)} \cos\frac{2nk\pi}{P} = \sum_{n=1}^P u_n \cos\frac{2nk\pi}{P} = \sum_{d\mid\gcd(k,P)} \mu(P_d)d \sum_{(a,b)\in Q_d} \cos\frac{2m_{ab}k\pi}{P}
\end{eqnarray}
\end{theorem}
\begin{proof}
First part of equation is given by definition of $u_n$. By taking $x = P$ in $E_1(x,s)$, we have
\begin{eqnarray} \label{E_1_P_s}
E_1(P,s) = \sum_{n = 1}^P u_n \cos 2ns = \sum_{d\mid P} \mu(P_d) \sum_{(a,b)\in Q_d} \frac{\sin B_{ab}(P)s - \sin A_{ab}s}{2\sin P_ds}
\end{eqnarray}
Since $B_{ab}(P) = A_{ab}+2P = [a,b]-2m_{ab}+2P$, then 
\begin{eqnarray}
\sin B_{ab}(P)s - \sin A_{ab}s & = & \cos A_{ab}s \sin 2Ps + \sin A_{ab}s(\cos 2Ps -1) \\
& = & \cos A_{ab}s \sin 2Ps - 2 \sin A_{ab}s \sin^2 Ps
\end{eqnarray}
Now we let $s\rightarrow \frac{k\pi}{P}$ in $E_1(P,s)$. It is clear that $\sin B_{ab}(P)s - \sin A_{ab}s\rightarrow 0$ as $s\rightarrow \frac{k\pi}{P}$. Let $d\mid P$. If $d\nmid k$, then $\sin P_ds\rightarrow \sin\frac{k\pi}{d} \ne 0$ as $s\rightarrow \frac{k\pi}{P}$. If $d\mid k$, then, by L'H\^{o}pital's rule, we have a limit: 
\begin{eqnarray}
\lim_{s\rightarrow\frac{k\pi}{P}} \frac{\sin B_{ab}(P)s - \sin A_{ab}s}{2\sin P_ds} & = & \lim_{s\rightarrow\frac{k\pi}{P}} \frac{\cos A_{ab}s \sin 2Ps}{2\sin P_ds} \\
& = & \left.\frac{2P\cos A_{ab}s\cos 2Ps}{2P_d \cos P_ds}\right|_{s = \frac{k\pi}{P}} \\
& = & d \cos\frac{k\pi}{d}\cos\frac{A_{ab}k\pi}{P} = d \cos\frac{2m_{ab}k\pi}{P} 
\end{eqnarray}
By taking $s\rightarrow \frac{k\pi}{P}$ in equation $(\ref{E_1_P_s})$, we have the theorem. 
\end{proof}

For each $p\mid P$, let $\bar p$ be the inverse of $p$ in $P/p = P_p$, that is the solution $z\in \mathcal{A}_{1,P}$ of systems: 
\begin{eqnarray}
zp \equiv 1 \mbox{ mod } P_p \mbox{ and } z \equiv 1 \mbox{ mod } p 
\end{eqnarray}
For $d\mid P$, a factor of $P$, let $\bar d = \prod_{p\mid d} \bar p$, then 
\begin{eqnarray}
\bar d d = \prod_{p\mid d} \bar pp \equiv 1 \mbox{ mod } P_d
\end{eqnarray}
and $\bar d\perp P$. For $d = 1$, we understand $\bar d = 1$. Let 
\begin{eqnarray}
\delta_2(a,b) = \left\{\begin{array}{ll}
1 & \mbox{if } 2\mid [a,b] \mbox{ and } N \mbox{ is odd}, \\ 0 & \mbox{otherwise}
\end{array}\right.
\end{eqnarray}
and $e(x) = e^{2\pi ix}$. 
\begin{theorem}
For $(a,b)\in Q$, let $a' = a/\gcd(2N,a)$, $b'=b/\gcd(2N,b)$ and  
\begin{eqnarray}
m_{ab}' = \delta_2(a,b)N\overline{P_2}P_2 + N \sum_{p\mid a'} \overline{P_p}P_p - N \sum_{p\mid b'} \overline{P_p}P_p
\end{eqnarray}
then $m_{ab}\equiv m_{ab}' \mbox{ \rm mod } [a,b]$. Let $k\geq 1$ and $d = P/[a,b]$. If $d\mid k$, then 
\begin{eqnarray}
e\left(\frac{m_{ab}k}{P}\right) = e\left(\delta_2(a,b)kN \frac{\overline{P_2}}{2} + kN\sum_{p\mid a'} \frac{\overline{P_p}}{p} - kN\sum_{p\mid b'}\frac{\overline{P_p}}{p}\right)
\end{eqnarray} 
\end{theorem}

\begin{proof}
First, we prove $a\mid N-m_{ab}'$. Assume prime $q\mid a$. If $q\mid N$, then $q\mid N-m_{ab}'$. Now assume $q\nmid N$. 
Notice that $q\mid P_p$ unless $p = q$ and $q\mid\overline{P_q}P_q-1$. If $q$ is odd, then $q\mid a'$ and $q\mid N-m_{ab}'$. If $q\nmid N$ and $q=2$, then $\delta_2(a,b) = 1$, $m_{ab}'$ is odd and $2\mid N-m_{ab}'$. Therefore, $a\mid N-m_{ab}'$. Similarly, we have $b\mid N+m_{ab}'$. Thus, $m_{ab}\equiv m_{ab}'$ mod $[a,b]$ and $m_{ab}k\equiv m_{ab}'k$ mod $P$. Since $P_p/P = 1/p$, then 
\begin{eqnarray}
e\left(\frac{m_{ab}k}{P}\right) = e\left(\frac{m_{ab}'k}{P}\right) = e\left(\delta_2(a,b)kN \frac{\overline{P_2}}{2} + kN\sum_{p\mid a'} \frac{\overline{P_p}}{p} - kN\sum_{p\mid b'}\frac{\overline{P_p}}{p}\right)
\end{eqnarray} 
This completes the proof.  
\end{proof}

\begin{theorem}
For $d\mid P$ and integer $k\geq 1$, if $d\mid k$, then 
\begin{eqnarray} \label{cos_mab_product}
\sum_{(a,b)\in Q_d}\cos\frac{2m_{ab}k\pi}{P} = (-1)^{k+kP_{dN}}\prod_{p\mid P_{2dN}} 2\cos\frac{2kN\overline{P_p}\pi}{p}
\end{eqnarray}
\end{theorem}
\begin{proof}
For $(a,b)\in Q_d$, let $a' = a/\gcd(2N,a)$ and $b' = b/\gcd(2N,b)$, then $P_{2dN} = [a',b']$. In fact, by the definition of $Q_d$, we have $b' = b$ and $P_{2dN} = a'b'$. If $P_{dN}$ is odd, then $P_{dN} = P_{2dN}$ and $\delta_2(a,b) = 0$. By the previous theorem, we have 
\begin{eqnarray}
\sum_{(a,b)\in Q_d} e\left(\frac{m_{ab}k}{P} \right) = \sum_{J(P_{2dN})} e\left(kN \sum_{p\mid P_{2dN}} j_p \frac{\overline{P_p}}{p} \right)
\end{eqnarray} 
where summation condition $J(P_{2dN})$ goes over all $j_p = \pm 1$ for each $p\mid P_{2dN}$. Thus, 
\begin{eqnarray}
\sum_{(a,b)\in Q_d} e\left(\frac{m_{ab}k}{P}\right) = \prod_{p\mid P_{2dN}} \left(e\left(kN \frac{\overline{P_p}}{p} \right) + e\left(-kN\frac{\overline{P_p}}{p} \right)\right) = \prod_{p\mid P_{2dN}} 2\cos\frac{2kN\overline{P_p}\pi}{p} 
\end{eqnarray} 
Since it is real, then 
\begin{eqnarray}
\sum_{(a,b)\in Q_d}\cos\frac{2m_{ab}k\pi}{P} = \sum_{(a,b)\in Q_d} e\left(\frac{m_{ab}k}{P}\right) = \prod_{p\mid P_{2dN}} 2\cos\frac{2kN\overline{P_p}\pi}{p} 
\end{eqnarray}
and the theorem is valid for odd $P_{dN}$. If $P_{dN}$ is even, then $N$ is odd, $\delta_2(a,b) = 1$ and 
\begin{eqnarray}
\sum_{(a,b)\in Q_d} e\left(\frac{m_{ab}k}{P} \right) & = & \sum_{J(P_{2dN})} e\left(\frac{kN\overline{P_2}}{2}+kN \sum_{p\mid P_{2dN}} j_p \frac{\overline{P_p}}{p} \right) \\
& = & e\left(\frac{kN\overline{P_2}}{2}\right) \sum_{J(P_{2dN})} e\left(kN \sum_{p\mid P_{2dN}} j_p \frac{\overline{P_p}}{p} \right)
\end{eqnarray} 
Since $e\left(\frac{kN\overline{P_2}}{2}\right) = (-1)^k$ when $P_{dN}$ is even, then we complete the proof. 
\end{proof}

\begin{theorem}
{\bf (Goldbach cosine sum-product formula on $W_P(N)$)} For $k\geq 1$, 
\begin{eqnarray}
\sum_{n\in W_P(N)} \cos\frac{2nk\pi}{P} = \mu(N_P)\prod_{p\mid\gcd(k,N_P)}(1-p) \prod_{p\mid P_{2N}} \alpha_p(P,kN)
\end{eqnarray}
\end{theorem}
where 
\begin{eqnarray}
\alpha_p(P,kN) = \left\{\begin{array}{ll}
p-2 & \mbox{if } p\mid k, \\ -2\cos\frac{2kN\overline{P_p}\pi}{p} & \mbox{if } p\nmid k
\end{array}\right.
\end{eqnarray}
\begin{proof}
By use of formula $(\ref{cos_mab_product})$ and $\mu(P_d) = \mu(P)\mu(d)$, equation $(\ref{cosine_k_sum_sieve_eq})$ becomes 
\begin{eqnarray}
\sum_{n\in W_P(N)} \cos\frac{2nk\pi}{P} = \mu(P)\sum_{d\mid\gcd(k,P)} (-1)^{k+kP_{dN}}\mu(d)d \prod_{p\mid P_{2dN}} 2\cos\frac{2kN\overline{P_p}\pi}{p}
\end{eqnarray} 
Since $\cos\frac{2kN\overline{P_p}\pi}{p} = 1$ for $p\mid k$, then    
\begin{eqnarray}
\sum_{n\in W_P(N)} \cos\frac{2nk\pi}{P} = \mu(P)\prod_{p\mid P_{2N}} 2\cos\frac{2kN\overline{P_p}\pi}{p} \sum_{d\mid \gcd(k, P)} z_k(d) 
\end{eqnarray} 
where  
\begin{eqnarray}
z_k(d) = (-1)^{k+kP_{dN}}\mu(d)d \prod_{p\mid \gcd(d,P_{2N})} \frac{1}{2}
\end{eqnarray}
Assume odd prime $p\mid \gcd(k,P)$. Thus, $(-1)^{k+kP_{dN}} = (-1)^{k+kP_{dpN}}$ for $d\mid \gcd(k,P_p)$. If $p\mid N_P$, then 
\begin{eqnarray}
\sum_{d\mid \gcd(k,P)} z_k(d) = \sum_{d\mid \gcd(k,P_p)} (z_k(d)+z_k(pd)) = (1-p) \sum_{d\mid \gcd(k,P_p)} z_k(d)
\end{eqnarray}
If $p\mid P_{2N}$, then 
\begin{eqnarray}
\sum_{d\mid \gcd(k,P)} z_k(d) = \sum_{d\mid \gcd(k,P_p)} (z_k(d)+z_k(pd)) = \left(1-\frac{p}{2}\right) \sum_{d\mid \gcd(k,P_p)} z_k(d)
\end{eqnarray}
Let 
\begin{eqnarray}
z_k' = \prod_{p\mid \gcd(k,N_P)\atop p \geq 3}(1-p) \sum_{d\mid \gcd(2,k,P)} (-1)^{k+kP_{dN}}\mu(d)d 
\end{eqnarray}
then 
\begin{eqnarray}
\sum_{n\in W_P(N)} \cos\frac{2nk\pi}{P} & = & \mu(P)z_k' \prod_{p\mid P_{2N}} 2\cos\frac{2kN\overline{P_p}\pi}{p} \prod_{p\mid \gcd(k,P_{2N})}\left(1 - \frac{p}{2}\right) \\
& = & \mu(P)z_k' \prod_{p\mid P_{2N}} -\alpha_p (P,kN) \\
& = & \mu(P)\mu(P_{2N}) z_k' \prod_{p\mid P_{2N}} \alpha_p (P,kN) 
\end{eqnarray}
It is easy to verify that 
\begin{eqnarray}
\sum_{d\mid \gcd(2,k,P)} (-1)^{k+kP_{dN}} \mu(d)d = \left\{\begin{array}{ll}
1 & \mbox{if } 2\nmid P, \\ (-1)^{k-1} & \mbox{if } 2\mid P \mbox{ and } 2\mid N, \\
-1 & \mbox{if } 2\mid P \mbox{ and } 2\nmid N 
\end{array}\right.
\end{eqnarray}
If $2\nmid \gcd(k,N_P)$, then we have 
\begin{eqnarray}
z_k' = (-1)^{P_N+1} \prod_{p\mid \gcd(k,N_P)}(1-p) 
\end{eqnarray}

If $2\mid \gcd(k,N_P)$, then $2\mid\gcd(k,N,P)$ and it is simple to verify the formula above is also valid. Since $\mu(P)\mu(P_{2N})(-1)^{P_N+1} = \mu(\gcd(P,2N))(-1)^{P_N+1} = \mu(N_P)$, we complete the proof. 
\end{proof}

Goldbach cosine sum-product formula has been verified numerically for $P \leq 30030 = 2\cdot 3\cdot 5\cdot 7 \cdot 11\cdot 13$. 
\begin{example}
Let $N = 4$ and $P = 15$, then $n = 3, 12$ and $15$ for $n\in W_P(N)$. Now we have $\mu(N_P) = 1$ and for $k \perp 15$, 
\begin{eqnarray}
\cos\frac{3\cdot 2k\pi}{15} + \cos\frac{12\cdot 2k\pi}{15} + \cos\frac{15\cdot 2k\pi}{15} = 4 \cos\frac{4k\pi}{3}\cos\frac{6k\pi}{5}
\end{eqnarray}
\end{example}

Since $\alpha_p(P,kN) = 1$ for any integer $k$ if $p = 3\mid P_{2N}$, then we also have
\begin{eqnarray}
C_P(N,k) \triangleq \sum_{n\in W_P(N)} \cos\frac{2nk\pi}{P} = \mu(N_P) \prod_{p\mid\gcd(k,N_P)}(1-p) \prod_{p\mid P_{6N}} \alpha_p(P,kN)
\end{eqnarray}

\section{Goldbach Decomposition Theorem}

Let us recall that $\bar d$ is the inverse of $d$ in $P_d$: $\bar dd \equiv 1$ mod $P_d$, and $V_P^d(N)$ is the set of $k\in P_{cd}^\perp$ such that $cdk\in W_P^d(N)$. Purpose of this section is to prove the Goldbach cosine sum-product formula over $W_P^d(N)$. 
\begin{theorem}
For $d\mid P_{6N}$, $V_P^d(N) = V_{P_d}^1(\bar dN)$. 
\end{theorem}

\begin{proof}
Let $c$ be the index of $(P,N)$. 
For $d\mid P_{6N}$, $k\in V_P^d(N)$ if and only if $(cdk - N)(cdk + N)\perp P$ and $k\in P_{cd}^\perp$. Since $(cdk - N)(cdk + N)\perp P$ is equivalent to $(cdk - N)(cdk + N)\perp P_d$ and is equivalent to $(ck - \bar dN)(ck + \bar dN)\perp P_d$, then $V_P^d(N) = V_{P_d}^1(\bar dN)$. 
\end{proof}

It is clear that $|W_P^d(N)| = |V_P^d(N)| = |V_{P_d}^1(\bar dN)| = |W_{P_d}^1(\bar dN)|$. 

\begin{lemma}
Assume $g(P,N)$ is defined for all square-free integers $P\geq 1$ and all integers $N\geq 1$. If 
\begin{eqnarray}
\sum_{d\mid P_{6N}} g(P_d,\bar dN) = 0
\end{eqnarray}
for each square-free integer $P\geq 1$ and each integer $N\geq 1$, then $g(P_d,\bar dN) = 0$ for any $d\mid P_{6N}$. 
\end{lemma}
\begin{proof}
We prove this lemma by induction on $\#P$, the number of prime factors in $P$. First, we take $P=1$ and $\#P=0$. In this case the lemma is obvious: $g(1,N) = 0$ for all integer $N\geq 1$. Now we assume the lemma is true for $\#P\leq n$, then $g(P_d,\bar dN) = 0$ if $\#P_d \leq n$ and $d\mid P_{6N}$. Let $P$ be square-free integer having $n+1$ prime factors. By assumption, we have
\begin{eqnarray}
\sum_{d\mid P_{6N}} g(P_d,\bar dN) = 0
\end{eqnarray}
and $g(P_d,\bar dN) = 0$ for any $d>1$ since $\#P_d \leq n$, then 
\begin{eqnarray}
0 = \sum_{d\mid P_{6N}} g(P_d,\bar dN) = g(P,N)
\end{eqnarray}
and the lemma is proved. 
\end{proof}

\begin{theorem}
{\bf (Goldbach decomposition theorem)} Assume $f(t)$ is defined for $0\leq t\leq 1$ and $g(P,N)$ is defined for all square-free integers $P\geq 1$ and all integers $N\geq 1$. If the following condition is satisfied: 
\begin{eqnarray}
\sum_{n\in W_P(N)} f\left(\frac{n}{P}\right) = \sum_{d\mid P_{6N}} g(P_d, \bar dN)
\end{eqnarray}
for each square-free integer $P\geq 1$ and each integer $N\geq 1$, then for any $d\mid P_{6N}$, 
\begin{eqnarray}
\sum_{n\in W_P^d(N)} f\left(\frac{n}{P}\right) = g(P_d, \bar dN)
\end{eqnarray}
\end{theorem}
\begin{proof}
Let 
\begin{eqnarray}
g'(P, N) = \sum_{n\in V_{P}^1(N)} f\left(\frac{n}{P_{c}}\right) = \sum_{n\in W_{P}^1(N)} f\left(\frac{n}{P}\right) 
\end{eqnarray}
then for $d\mid P_{6N}$, we have
\begin{eqnarray}
g'(P_d, \bar dN) = \sum_{n\in V_{P_d}^1(\bar dN)} f\left(\frac{n}{P_{cd}}\right) = \sum_{n\in V_P^d(N)} f\left(\frac{n}{P_{cd}}\right) = \sum_{n\in W_{P}^d(N)} f\left(\frac{n}{P}\right) 
\end{eqnarray}
and 
\begin{eqnarray}
\sum_{n\in W_P(N)} f\left(\frac{n}{P}\right) = \sum_{d\mid P_{6N}} \sum_{n\in W_P^d(N)} f\left(\frac{n}{P}\right) = \sum_{d\mid P_{6N}} g'(P_d, \bar dN) 
\end{eqnarray}
By the assumption of this theorem, we have
\begin{eqnarray}
\sum_{n\in W_P(N)} f\left(\frac{n}{P}\right) = \sum_{d\mid P_{6N}} g(P_d, \bar dN) = \sum_{d\mid P_{6N}} g'(P_d, \bar dN)
\end{eqnarray}
Thus, by the previous lemma, we have $g(P_d,\bar dN) = g'(P_d,\bar dN)$ for all $d\mid P_{6N}$. 
\end{proof}

\begin{lemma}
For $d\mid P_{6N}$, $p\mid P_{6dN}$ and $k\geq 1$, $\alpha_p(P,kN) = \alpha_p(P_d,k\bar dN)$. 
\end{lemma}
\begin{proof}
If $p\mid k$, then $\alpha_p(P,kN) = \alpha_p(P_d,k\bar dN) = p-2$. Now we assume $p\nmid k$. Let $P' = P_d$ and $P_p' = P' / p$, then $P_p' = P_{dp}$. Let $\overline{P_p'}$ be the inverse of $P_p'$ in $p$: $P_p'\overline{P_p'} \equiv 1$ mod $p$. Since
\begin{eqnarray}
\overline{P_p} = \overline{dP_{dp}} = \bar d\cdot \overline{P_{dp}} \equiv \bar d\cdot\overline{P_{dp}}P_p'\overline{P_p'} = \bar d\cdot\overline{P_{dp}}P_{dp}\cdot\overline{P_p'} \equiv \bar d\cdot\overline{P_p'} \mbox{ mod } p
\end{eqnarray}
then by the definition, we have
\begin{eqnarray}
\alpha_p(P,kN) = -2\cos\frac{kN\overline{P_p}\pi}{p} = -2\cos\frac{k\bar dN\overline{P_p'}\pi}{p} = \alpha_p(P_d,k\bar dN)
\end{eqnarray}
Thus, the proof is complete. 
\end{proof}

\begin{theorem} 
{\bf (Goldbach cosine sum-product formula on $W_P^d(N)$)} For $d\mid P_{6N}$ and $k\geq 1$, 
\begin{eqnarray}
C_P^d(N, k) \triangleq \sum_{n\in W_P^d(N)} \cos\frac{2nk\pi}{P} = \mu(N_P) \prod_{p\mid\gcd(k,N_P)}(1-p) \prod_{p\mid P_{6dN}} (\alpha_p(P,kN)-1)
\end{eqnarray}
\end{theorem}
\begin{proof}
Since 
\begin{eqnarray}
\sum_{n\in W_P(N)} \cos\frac{2nk\pi}{P} & = & \mu(N_P)\prod_{p\mid\gcd(k,N_P)}(1-p) \prod_{p\mid P_{6N}} \alpha_p(P,kN) \\
& = & \mu(N_P)\prod_{p\mid\gcd(k,N_P)}(1-p) \sum_{d\mid P_{6N}} \prod_{p\mid P_{6dN}} (\alpha_p(P,kN) - 1)
\end{eqnarray}
Let 
\begin{eqnarray}
h(P,N,k) = \mu(N_P) \prod_{p\mid\gcd(k,N_P)}(1-p) 
\end{eqnarray}
then $h(P_d,\bar dN,k) = h(P,N,k) \triangleq a_k$ is constant for $d\mid P_{6N}$. Since $\alpha_p(P,kN) = \alpha_p(P_d,k\bar dN)$ for any $d\mid P_{6pN}$, then 
\begin{eqnarray}
g(P_d, \bar dN) \triangleq a_k\prod_{p\mid P_{6d\bar dN}} \left(\alpha_p(P_d,k\bar dN) - 1\right) = a_k \prod_{p\mid P_{6dN}} (\alpha_p(P,kN)-1)
\end{eqnarray} 
and  
\begin{eqnarray}
\sum_{n\in W_P(N)} \cos\frac{2nk\pi}{P} = \sum_{d\mid P_{6N}} g(P_d, \bar dN)
\end{eqnarray}
By the Goldbach decomposition theorem, we have
\begin{eqnarray}
\sum_{n\in W_P^d(N)} \cos\frac{2nk\pi}{P} = g(P_d, \bar dN) = \mu(N_P) \prod_{p\mid\gcd(k,N_P)}(1-p) \prod_{p\mid P_{6dN}} (\alpha_p(P,kN)-1)
\end{eqnarray}
This completes the proof. 
\end{proof}

\begin{theorem}
For $d\mid P_{6N}$ and $k\geq 1$, $C_P^d(N, k) = C_{P_d}^1(\bar dN, k)$.
\end{theorem}
This is because
\begin{eqnarray}
\sum_{n\in W_P^d(N)} \cos\frac{2nk\pi}{P} = \sum_{n\in V_P^d(N)} \cos\frac{2nk\pi}{P_{cd}} = \sum_{n\in V_{P_d}^1(\bar dN)} \cos\frac{2nk\pi}{P_{cd}} = \sum_{n\in W_{P_d}^1(\bar dN)} \cos\frac{2nk\pi}{P_{d}}
\end{eqnarray}

\section{Goldbach Counting Function}

For $x >0$, let $S_P(N,x)$ be the number of $n\in \mathcal{I}_{1,x}$ such that $(N-n)(N+n)\perp P$. For $d\mid P_{6N}$, $S_P^d(N,x)$ be the number of $n\in \mathcal{I}_{1,x}$ such that $(N-n)(N+n)\perp P$ and $\gcd(P_{6N},n) = d$. Both of them are called the Goldbach counting functions. We define 
\begin{eqnarray}
t_P(N) = \left\{\begin{array}{cl} 
\frac{1}{2} & \mbox{if } P\in W_P(N), \\ 0 & \mbox{if } P \not\in W_P(N)
\end{array}\right. \mbox{ and } T_P(N,x) = \sum_{n\in W_P(N)} \left(\left\{\frac{x-n}{P}\right\} -\frac{1}{2}\right)
\end{eqnarray}

\begin{theorem}
Goldbach counting function can be given as follows: 
\begin{eqnarray}
S_P(N,x) = |W_P(N)|\frac{x}{P} - t_P(N) - T_P(N,x) 
\end{eqnarray}
\end{theorem}

\begin{proof}
We only need to prove the formula for $0 < x \leq P$. For $0 < x \leq P$, we have 
\begin{eqnarray} 
S_P(N,x) & = & \sum_{n\in W_P(N)} \left\lfloor\frac{P+x-n}{P}\right\rfloor \\
& = & \sum_{n\in W_P(N)} \frac{P+x-n}{P} - \sum_{n\in W_P(N)} \left\{\frac{P+x-n}{P}\right\} \\
& = & |W_P(N)|\frac{x}{P} + \sum_{n\in W_P(N)} \frac{P-n}{P} - \sum_{n\in W_P(N)} \left\{\frac{x-n}{P}\right\} 
\end{eqnarray}
Since $P-n\in W_P(N)$ if $1\leq n < P$, then 
\begin{eqnarray}
\sum_{n\in W_P(N)} \frac{P-n}{P} = \sum_{n\in W_P(N)} \frac{P-n}{2P} + \sum_{n\in W_P(N)} \frac{n}{2P} - t_P(N) = \sum_{n\in W_P(N)} \frac{1}{2} - t_P(N) 
\end{eqnarray}
That completes the proof. 
\end{proof}

By this theorem, we now extend the range of $x$ in $T_P(N,x)$ and $S_P(N,x)$ to the whole reals. 
\begin{theorem}
If $Pt$ is not integer, then 
\begin{eqnarray}
T_P(N,Pt) = - \sum_{k = 1}^\infty C_P(N,k) \frac{\sin 2k\pi t}{k\pi} 
\end{eqnarray}
\end{theorem}
\begin{proof}
For $0<t<1$, we have
\begin{eqnarray}
\{-t\} -\frac{1}{2} = 1 - t - \frac{1}{2} = - \left(\{t\} - \frac{1}{2}\right)
\end{eqnarray}
Thus, $\{t\} - \frac{1}{2}$ is an odd function of $t$ and its Fourier transform has only sine terms. We calculate its coefficient as follows for $k\geq 1$: 
\begin{eqnarray}
2\int_0^1 \left(t - \frac{1}{2}\right) \sin (2k\pi t)dt
= -\left.\frac{1}{k\pi} t\cos(2k\pi t)\right|_{t=0}^{t=1} - \frac{1}{k\pi} \int_0^1 \cos(2k\pi t)dt = -\frac{1}{k\pi}
\end{eqnarray}

Thus, if $t \ne \frac{n}{P}$ for $n\in W_P(N)$, or more stronger condition: if $Pt$ is not integer, then 
\begin{eqnarray} 
T_P(N,Pt) = \sum_{n\in W_P(N)} \left(\left\{t-\frac{n}{P}\right\} -\frac{1}{2}\right) = -\sum_{k = 1}^\infty \frac{1}{k\pi}\sum_{n\in W_P(N)} \sin 2k\left(t - \frac{n}{P}\right)\pi
\end{eqnarray}

Since 
\begin{eqnarray}
\sin 2k\left(t - \frac{n}{P}\right)\pi = \cos \frac{2nk\pi}{P} \sin 2k\pi t - \sin \frac{2nk\pi}{P} \cos 2k\pi t
\end{eqnarray}
and $P-n\in W_P(N)$ if $n\in W_P(N)$ and $n>0$, then 
\begin{eqnarray}
T_P(N,Pt) = -\sum_{k = 1}^\infty \frac{\sin 2k\pi t}{k\pi} \sum_{n\in W_P(N)} \cos \frac{2nk\pi}{P} = - \sum_{k = 1}^\infty C_P(N,k) \frac{\sin 2k\pi t}{k\pi} 
\end{eqnarray}
Since $T_P(N,Pt)$ is a periodic function of $t$ with period $1$, then we have the theorem. 
\end{proof}

For $d\mid P_{6N}$, we define 
\begin{eqnarray} 
T_P^d(N,Pt) = \sum_{n\in W_P^d(N)} \left(\left\{t - \frac{n}{P}\right\} -\frac{1}{2}\right)
\end{eqnarray}
By calculation of its Fourier coefficients, if $Pt$ is not integer, we have 
\begin{eqnarray} 
T_P^d(N,Pt) = - \sum_{k = 1}^\infty C_P^d(N,k) \frac{\sin 2k\pi t}{k\pi} 
\end{eqnarray}

Notice that $P\in W_P(N)$ if and only if $N_P = 1$. Let 
\begin{eqnarray}
\delta_6(n,m) = \left\{\begin{array}{cl}
\frac{1}{2} & \mbox{if } n\mid 6 \mbox{ and } n\perp m, \\ 0 & \mbox{otherwise} 
\end{array}\right. 
\end{eqnarray}

and $t_P^d(N) = \delta_6(P_d,\bar dN)$ for $d\mid P_{6N}$. 
\begin{lemma} {\bf (Decomposition of $t_P(N)$)}
\begin{eqnarray}
t_P(N) = \sum_{d\mid P_{6N}} t_P^d(N) = \sum_{d\mid P_{6N}} \delta_6(P_d,\bar dN)
\end{eqnarray}
\end{lemma}
\begin{proof}
Let $d' = P_{6N}$. If $d\mid P_{6N}$ and $d < P_{6N}$, then $P_d \nmid 6$, $\delta_6(P_d,\bar dN) = 0$ and 
\begin{eqnarray}
\sum_{d\mid P_{6N}} \delta_6(P_d,\bar dN) = \delta_6(P_{d'},\bar d'N)
\end{eqnarray}
Since $P_{d'} = P/d'=P/P_{6N} = \gcd(P,6N)$ and $\bar d'd'\perp P_{d'}$, then 
\begin{eqnarray}
\gcd(P_{d'},\bar d'N) = \gcd(P_{d'},N) = \gcd(P,6N,N) = \gcd(P,N) = N_P 
\end{eqnarray}
If $t_P(N) = 0$, then $N_P > 1$ and $\delta_6(P_{d'},\bar d'N)=0$. If $t_P(N) = \frac{1}{2}$, then $N_P = 1$. Thus, $P_{d'} = \gcd (P,6)\mid 6$ and $P_{d'}\perp \bar dN$. By definition, $\delta_6(P_{d'},\bar d'N)=\frac{1}{2}$. Hence, $\delta_6(P_{d'},\bar d'N)=t_P(N)$. 
\end{proof}

\begin{theorem}
For $d\mid P_{6N}$, $S_P^d(N,Pt) = S_{P_d}^1(\bar dN,P_dt)$ and $T_P^d(N,Pt) = T_{P_d}^1(\bar dN,P_dt)$. 
\end{theorem}
\begin{proof}
We need to prove the theorem only for $0 < t < 1$ and now we assume $0 < t < 1$. Since $cdk \in W_P^d(N)$ and $cdk \leq Pt$ if and only if $ck\in W_{P_d}^1(\bar dN)$ and $ck\leq P_dt$, then $S_P^d(N,Pt) = S_{P_d}^1(\bar dN,P_dt)$. Since $C_P^d(N,k) = C_{P_d}^1(\bar dN,k)$ for $d\mid P_{6N}$, then $T_P^d(N,Pt) = T_{P_d}^1(\bar dN,P_dt)$. 
\end{proof}

\begin{corollary}
For $d\mid P_{6N}$ and $d'\mid P_{6dN}$, $T_{P_d}^{d'}(\bar dN,P_dt) = T_{P_{dd'}}^1(\overline{dd'}N,P_{dd'}t) = T_P^{dd'}(N,Pt)$. 
\end{corollary}

\begin{theorem} {\bf (Decomposition of $S_P(N,Pt)$)} For $d\mid P_{6N}$, 
\begin{eqnarray} \label{decomp_d_p_n_x}
S_P^d(N,Pt) = |W_P^d(N)|t - t_P^d(N) - T_P^d(N,Pt) 
\end{eqnarray}
\end{theorem}
\begin{proof}
Since we have the following: 
\begin{eqnarray}
S_P(N,Pt) & = & |W_P(N)| t - t_P(N) - T_P(N,Pt) \\
\sum_{d\mid P_{6N}} S_P^d(N,Pt) & = & \sum_{d\mid P_{6N}} \left(|W_P^d(N)|t - t_P^d(N) - T_P^d(N,Pt)\right) \\
\sum_{d\mid P_{6N}} S_{P_d}^1(\bar dN,P_dt) & = & \sum_{d\mid P_{6N}} \left(|W_{P_d}^1(\bar dN)|t - \delta_6(P_d,\bar dN)  - T_{P_d}^1(\bar dN,P_dt)\right)
\end{eqnarray}
Then the theorem follows by Goldbach decomposition theorem. 
\end{proof}

By this theorem, again, for each $d\mid P_{6N}$, we now extend the range of $x$ in $T_P^d(N,x)$ and $S_P^d(N,x)$ to the whole reals.

\section{Deduction Formula for Goldbach Counting Function}

Let us start with deduction formula for $T_P(N,x)$ in terms of $T_{P_p}(N,x)$ for $p\mid P_{2N}$. 
\begin{lemma}
$C_P(N,k) = C_{P_p}(\bar pN,k)\alpha_p(P,kN)$ for $p\mid P_{2N}$ and $k\geq 1$. 
\end{lemma}
\begin{proof}
Let $P' = P_p$ and $P_n' = P'/\gcd(P',n)$. Thus, we have
\begin{eqnarray}
C_P(N,k) & = & \alpha_p(P,kN)\mu(N_P) \prod_{q\mid\gcd(k,N_P)}(1-q) \prod_{q\mid P_{2pN}} \alpha_q(P,kN) \\
& = & \alpha_p(P,kN)\mu(N_{P'}) \prod_{q\mid\gcd(k,N_{P'})}(1-q) \prod_{q\mid P_{2N}'} \alpha_q(P,kN) 
\end{eqnarray} 
Let $\overline{P_q'}$ be the inverse of $P_q'$ in $q$ for $q\mid P_{2N}' = P_{2pN}$: $\overline{P_q'} P_q' \equiv 1$ mod $q$. Since
\begin{eqnarray}
\bar p\overline{P_q'} \equiv \bar p\overline{P_q'} P_{pq} \overline{P_{pq}} = \bar p\overline{P_q'} P_q' \overline{P_{pq}} \equiv \bar p\cdot\overline{P_{pq}} = \overline{pP_{pq}} = \overline{P_q} \mbox{ mod } q
\end{eqnarray}
then $\alpha_q(P,kN) = \alpha_q(P',k\bar pN)$. This completes the proof. 
\end{proof}
Similarly, we have
\begin{lemma}
$C_P^d(N,k) = C_{P_p}^d(\bar pN,k)(\alpha_p(P,kN)-1)$ for $d\mid P_{6N}$, $p\mid P_{6dN}$ and $k\geq 1$. 
\end{lemma}
\begin{theorem}  {\bf (First deduction formula for $T_P(N,x)$)} 
If $p\mid P_{2N}$ and $x$ is not integer, then
\begin{eqnarray}
T_P(N,x) = T_{P_p}(N,x) - T_{P_p}\left(\bar pN, \frac{N\overline{P_p}P_p+x}{p}\right) + T_{P_p}\left(\bar pN, \frac{N\overline{P_p}P_p-x}{p}\right)
\end{eqnarray}
\end{theorem}
\begin{proof}
For $p\mid P_{2N}$, we have $C_P(N,k) = C_{P_p}(\bar pN,k)\alpha_p(P,kN)$. Let $t = x/P$ and $y = N\overline{P_p}/p$. Since $x$ is not integer, then by the definition of $\alpha_p(P,kN)$, 
\begin{eqnarray}
T_P(N,x) & = & -\sum_{k=1}^\infty C_{P_p}(\bar pN,k)\alpha_p(P,kN) \frac{\sin 2k\pi t}{k\pi} \\
& = & 2\sum_{k=1\atop p\nmid k}^\infty C_{P_p}(\bar pN,k) \cos 2ky\pi \frac{\sin 2k\pi t}{k\pi} - \sum_{k=1}^\infty C_{P_p}(\bar pN,pk)(p-2)\frac{\sin 2pk\pi t}{pk\pi}\\ 
& = & 2\sum_{k=1}^\infty C_{P_p}(\bar pN,k) \cos 2ky\pi \frac{\sin 2k\pi t}{k\pi} - \sum_{k=1}^\infty C_{P_p}(\bar pN,pk) \frac{\sin 2pk\pi t}{k\pi}
\end{eqnarray}
Since $\bar pp \equiv 1$ mod $q$ for any prime $q\mid P_{2pN}$, then $C_{P_p}(\bar pN,pk) = C_{P_p}(N,k)$ and 
\begin{eqnarray}
T_P(N,x) & = & 2\sum_{k=1}^\infty C_{P_p}(\bar pN,k) \cos 2ky\pi \frac{\sin 2k\pi t}{k\pi} - \sum_{k=1}^\infty C_{P_p}(N,k)\frac{\sin 2pk\pi t}{k\pi} \\
& = & 2\sum_{k=1}^\infty C_{P_p}(\bar pN,k) \cos 2ky\pi \frac{\sin 2k\pi t}{k\pi} + T_{P_p}(N,x)
\end{eqnarray}

Since $2\cos 2ky\pi \sin 2k\pi t = \sin 2k\pi(y+t) - \sin 2k\pi(y-t)$, 
\begin{eqnarray}
P_p(y+t) = \frac{N\overline{P_p}P_p+x}{p} \mbox{ and } P_p(y-t) = \frac{N\overline{P_p}P_p-x}{p}
\end{eqnarray}
then we have the proof. 
\end{proof}

\begin{theorem} {\bf (First deduction formula for $S_P(N,x)$)} 
If $p\mid P_{2N}$ and $x$ is not integer, then 
\begin{eqnarray}
S_P(N,x) = S_{P_p}(N,x) - S_{P_p}\left(\bar pN, \frac{N\overline{P_p}P_p+x}{p}\right) + S_{P_p}\left(\bar pN, \frac{N\overline{P_p}P_p-x}{p}\right)
\end{eqnarray}
\end{theorem}
\begin{proof}
Let $t = x/P$ and $z = N\overline{P_p}P_p$. Since $|W_P(N)| = |W_{P_p}(N)|(p-2)$ and $|W_{P_p}(\bar pN)| = |W_{P_p}(N)|$, then
\begin{eqnarray}
S_P(N,x) & = & |W_{P_p}(N)|(p-2) t - t_P(N) - T_P(N,x) \\
S_{P_p}(N,x) & = & |W_{P_p}(N)|pt - t_{P_p}(N) - T_{P_p}(N,x) \\
S_{P_p}\left(\bar pN,\frac{z+x}{p}\right) & = & |W_{P_p}(N)| \frac{z+x}{P} - t_{P_p}(\bar pN) - T_{P_d}\left(\bar pN,\frac{z+x}{p}\right) \\
S_{P_p}\left(\bar pN,\frac{z-x}{p}\right) & = & |W_{P_p}(N)| \frac{z-x}{P} - t_{P_p}(\bar pN) - T_{P_d}\left(\bar pN,\frac{z-x}{p}\right)
\end{eqnarray}
Since $t_P(N) = t_{P_p}(N)$, then by the previous theorem, we have the formula. 
\end{proof}
Again, this formula is verified numerically for several cases. 
\begin{theorem} {\bf (First deduction formula for $T_P^d(N,x)$)}
If $p\mid P_{6N}$, $d\mid P_{6pN}$ and $x$ is not integer, then 
\begin{eqnarray}
T_P^d(N,x) = T_{P_p}^d(N,x) - T_P^{dp}(N,x) - T_P^{dp}(N,N\overline{P_p}P_p+x) + T_P^{dp}(N,N\overline{P_p}P_p-x)
\end{eqnarray}
\end{theorem}
\begin{proof}
For $d\mid P_{6pN}$, $C_P^d(N,k) = C_{P_p}^d(\bar pN,k)(\alpha_p(P,kN)-1)$. Let $t = x/P$, then 
\begin{eqnarray}
T_P^d(N,x) & = & -\sum_{k=1}^\infty C_{P_p}^d(\bar pN,k)(\alpha_p(P,kN) - 1) \frac{\sin 2k\pi t}{k\pi} \\
& = & -\sum_{k=1}^\infty C_{P_p}^d(\bar pN,k)\alpha_p(P,kN) \frac{\sin 2k\pi t}{k\pi} - T_{P_p}^d(\bar pN,P_pt)
\end{eqnarray} 
Similar to first deduction formula for $T_P(N,x)$, let $y = N\overline{P_p}/p$, then we have
\begin{eqnarray}
 &  & -\sum_{k=1}^\infty C_{P_p}^d(\bar pN,k)\alpha_p(P,kN) \frac{\sin 2k\pi t}{k\pi} \\
& = & 2\sum_{k=1\atop p\nmid k}^\infty C_{P_p}^d(\bar pN,k) \cos 2ky\pi \frac{\sin 2k\pi t}{k\pi} - \sum_{k=1}^\infty C_{P_p}^d(\bar pN,pk)(p-2)\frac{\sin 2pk\pi t}{pk\pi}\\ 
& = & 2\sum_{k=1}^\infty C_{P_p}^d(\bar pN,k) \cos 2ky\pi \frac{\sin 2k\pi t}{k\pi} - \sum_{k=1}^\infty C_{P_p}^d(\bar pN,pk) \frac{\sin 2pk\pi t}{k\pi}
\end{eqnarray}
Since $2\cos 2ky\pi \sin 2k\pi t = \sin 2k\pi(y+t) - \sin 2k\pi(y-t)$ and $C_{P_p}^d(\bar pN,pk) = C_{P_p}^d(N,k)$, then
\begin{eqnarray}
T_P^d(N,x) & = & T_{P_p}^d(N,x) - T_{P_p}^d\left(\bar pN, \frac{N\overline{P_p}P_p+x}{p}\right) + T_{P_p}^d\left(\bar pN, \frac{N\overline{P_p}P_p-x}{p}\right) - T_{P_p}^d\left(\bar pN,\frac{x}{p}\right) \\
& = & T_{P_p}^d(N,x) - T_P^{dp}(N,N\overline{P_p}P_p+x) + T_P^{dp}(N,N\overline{P_p}P_p-x) - T_P^{dp}(N,x)
\end{eqnarray}
This completes the proof. 
\end{proof}

We have a similar result for $p\mid N_P$ as follows: 
\begin{theorem} {\bf (Second deduction formula for $T_P(N,x)$)} 
If $p\mid N_P$ and $x$ is not integer, then 
\begin{eqnarray}
T_P(N,x) = T_{P_p}(N,x) - T_{P_p}\left(\bar pN,\frac{x}{p}\right)
\end{eqnarray}
\end{theorem}
\begin{proof}

If $p\mid N_P$, then $\mu(N_P) = -\mu(N_{P_p})$ and 
\begin{eqnarray}
C_P(N,k) & = & -\mu(N_{P_p}) \prod_{q\mid\gcd(k,N_P)}(1-q) \prod_{q\mid P_{6N}} \alpha_q(P,kN) \\
& = & -\mu(N_{P_p}) p_k \prod_{q\mid\gcd(k,N_{P_p})}(1-q) \prod_{q\mid P_{6p\bar pN}} \alpha_q(P_p,k\bar pN) \\
& = & -C_{P_p}(\bar pN,k)p_k
\end{eqnarray}
where $p_k = 1$ if $p\nmid k$ and $p_k = 1-p$ if $p\mid k$. Thus, 
\begin{eqnarray}
T_P(N,x) & = & \sum_{k=1}^\infty C_{P_p}(\bar pN,k)p_k \frac{\sin 2k\pi t}{k\pi} \\
& = & \sum_{k=1\atop p\nmid k}^\infty C_{P_p}(\bar pN,k) \frac{\sin 2k\pi t}{k\pi} + (1-p)\sum_{k=1}^\infty C_{P_p}(\bar pN,pk)\frac{\sin 2pk\pi t}{pk\pi} \\
& = & \sum_{k=1}^\infty C_{P_p}(\bar pN,k) \frac{\sin 2k\pi t}{k\pi} - \sum_{k=1}^\infty C_{P_p}(\bar pN,pk)\frac{\sin 2pk\pi t}{k\pi} \\ 
& = & -T_{P_p}\left(\bar pN,\frac{x}{p}\right) + T_{P_p}(N,x)
\end{eqnarray}
and we have the proof. 
\end{proof}

Since $C_P^d(N,k) = -C_{P_p}^d(\bar pN,k)p_k$ for $p\mid N_P$ and $d\mid P_{6N}$, then 
\begin{theorem} {\bf (Second deduction formula for $T_P^d(N,x)$)} 
If $p\mid N_P$ and $x$ is not integer, then for $d\mid P_{6N}$, 
\begin{eqnarray}
T_P^d(N,x) = T_{P_p}^d(N,x) - T_{P_p}^d\left(\bar pN,\frac{x}{p}\right)
\end{eqnarray}
\end{theorem}

\begin{theorem} {\bf (Third deduction formula for $T_P(N,x)$)} 
Let $c$ be the index of $(P,N)$. If $x$ is not integer and $d\mid P_{6N}$, then 
\begin{eqnarray}
T_P(N,x) = T_{P_c}\left(\bar cN,\frac{x}{c}\right) \mbox{ and } T_P^d(N,x) = T_{P_c}^d\left(\bar cN,\frac{x}{c}\right)
\end{eqnarray}
\end{theorem}

This is because $C_P(N,k) = C_{P_c}(\bar cN,k)$ and $C_P^d(N,k) = C_{P_c}^d(\bar cN,k)$ for any integer $k$.

\section{Densities of $W_P(N)$ and $W_P^d(N)$}

We give an estimate of $|W_P(N)|/P$ when $P$ is the product of all primes $\leq z$. First let us define 
\begin{eqnarray}
C_1 = e^\gamma \approx 1.781072418
\end{eqnarray}

where $\gamma \approx 0.5772156649$ is Euler-Mascheroni constant. Let 
\begin{eqnarray}
C_2 = \prod_{p\geq 3}\frac{(p-2)p}{(p-1)^2} \approx 0.66016 \mbox{ and } d_N = \prod_{p\mid N \atop p\geq 3}\frac{p-1}{p-2}
\end{eqnarray}

Let $\omega_P(N) = |W_P(N)|/P$, the density of $W_P(N)$ between 1 and $P$, then  
\begin{eqnarray}
\omega_P(N) = \frac{|W_P(N)|}{P} = \frac{1}{c'}\prod_{p\mid N_P}\left(1-\frac{1}{p}\right) \prod_{p\mid P_{2N}} \left(1-\frac{2}{p}\right)
\end{eqnarray}
where $c' = \frac{\gcd(NP,2)}{\gcd(N,2)}$.  
\begin{theorem}
If $P$ is the product of all primes $\leq z$, then 
\begin{eqnarray}
\omega_P(N) = 2 \prod_{3\leq p\leq t} \frac{(p-2)p}{(p-1)^2} \prod_{p\mid N\atop 3\leq p \leq z} \frac{p-1}{p-2} \prod_{2\leq p\leq z} \left(1-\frac{1}{p}\right)^2 \sim \frac{2C_2d_N}{C_1^2\log^2 z}
\end{eqnarray}
The asymptotic formula holds for large $z$.  
\end{theorem}
\begin{proof}
By definition, we have
\begin{eqnarray}
\omega_P(N) = \frac{1}{2}\prod_{p\mid N\atop 3\leq p \leq z} \frac{p-1}{p-2} \prod_{3\leq p\leq z}\left(1-\frac{2}{p}\right)
\end{eqnarray}
By the following identities 
\begin{eqnarray}
\left(1-\frac{2}{p}\right) & = & \left(1-\frac{1}{p}\right)^2 \frac{(p-2)p}{(p-1)^2} \\
\prod_{3\leq p\leq z}\left(1-\frac{2}{p}\right) & = & 4\prod_{2\leq p\leq z} \left(1-\frac{1}{p}\right)^2 \prod_{3\leq p\leq z} \frac{(p-2)p}{(p-1)^2} 
\end{eqnarray}
we have 
\begin{eqnarray}
\omega_P(N) = 2 \prod_{3\leq p\leq z} \frac{(p-2)p}{(p-1)^2} \prod_{p\mid N\atop 3\leq p \leq z} \frac{p-1}{p-2} \prod_{2\leq p\leq z} \left(1-\frac{1}{p}\right)^2
\end{eqnarray}

First product approaches to $C_2$ as $z$ approaches to infinity and second product approaches to $d_N$. Third product approaches to $e^{-2\gamma}\log^{-2} z$ as the Mertens' 3rd theorem gives 
\begin{eqnarray}
\lim_{y\rightarrow \infty} \log y \prod_{2\leq p\leq y} \left(1-\frac{1}{p}\right) = e^{-\gamma}
\end{eqnarray}
That completes the proof.
\end{proof}

Hardy and Littlewood made the following assertion: asymptotic formula \cite{Hardy} and \cite{Wang}: 
\begin{eqnarray}
S_P(N,N) \sim \frac{2C_2d_N}{\log^2 N}N 
\end{eqnarray}

where $z = \sqrt{2N}$. Since $\omega_P(N)N = S_P(N,N)+T_P(N,N)+t_P(N)$ and 
\begin{eqnarray}
\omega_P(N)N \sim \frac{8C_2d_N}{C_1^2\log^2 N}N 
\end{eqnarray}

then, we will have, if Hardy and Littlewood are correct, 
\begin{eqnarray}
\frac{T_P(N,N)}{S_P(N,N)} \approx \frac{\omega_P(N)N - S_P(N,N)}{S_P(N,N)} \sim \frac{8/C_1^2-2}{2} \approx 0.260947
\end{eqnarray}

It means the error term $T_P(N,N)$ holds 26\% of target function $S_P(N,N)$. Numerical calculation strongly supports this assertion. \\ \\
For $d\mid P_{6N}$, the density of $W_P^d(N)$ between $1$ and $P$ is defined as 
\begin{eqnarray}
\omega_P^d(N) \triangleq \frac{|W_P^d(N)|}{P} = \frac{1}{cd} \prod_{p\mid N_P}\left(1-\frac{1}{p}\right) \prod_{p\mid P_{6dN}} \left(1-\frac{3}{p}\right)
\end{eqnarray}

where $c = \frac{\gcd(NP,6)}{\gcd(N,6)}$, the index of $(P,N)$. By use of this density, we have
\begin{theorem} 
For $d\mid P_{6N}$, $S_P^d(N,x) = x\omega_P^d(N) - t_P^d(N) - T_P^d(N,x)$. 
\end{theorem}
Let 
\begin{eqnarray}
C_3 = \prod_{p\geq 5} \frac{(p-3)p^2}{(p-1)^3} \approx 0.635166 \mbox{ and } d_N' = \prod_{3\mid N} 2 \prod_{p\mid N\atop p \geq 5} \frac{p-1}{p-3} 
\end{eqnarray}

\begin{theorem}
If $P$ is the product of all primes $\leq z$, then density 
\begin{eqnarray}
\omega_P^1(N) = \frac{9}{2} \prod_{5\leq p\leq z} \frac{(p-3)p^2}{(p-1)^3} \prod_{3\mid N} 2 \prod_{p\mid N\atop 5\leq p \leq z} \frac{p-1}{p-3} \prod_{2\leq p\leq z}\left(1-\frac{1}{p}\right)^3 \sim \frac{9C_3d_N'}{2C_1^3\log^3 z} 
\end{eqnarray} 
The asymptotic formula holds for large $z$. 
\end{theorem}
\begin{proof}
Now we have $\gcd(NP,6)=6$ and
\begin{eqnarray}
\omega_P^1(N) & = & \frac{\gcd(N,6)}{6} \prod_{p\mid \gcd(N,6)} \frac{p-1}{p}\prod_{p\mid N\atop 5\leq p \leq z} \frac{p-1}{p-3} \prod_{5\leq p\leq z}\left(1-\frac{3}{p}\right) \\
& = & \frac{1}{6} \prod_{3\mid N} 2 \prod_{p\mid N\atop 5\leq p \leq z} \frac{p-1}{p-3} \prod_{5\leq p\leq z}\left(1-\frac{3}{p}\right) 
\end{eqnarray}
By the following identities 
\begin{eqnarray}
\left(1-\frac{3}{p}\right) & = & \left(1-\frac{1}{p}\right)^3 \frac{(p-3)p^2}{(p-1)^3} \\
\prod_{5\leq p\leq z}\left(1-\frac{3}{p}\right) & = & 27 \prod_{2\leq p\leq z} \left(1-\frac{1}{p}\right)^3 \prod_{5\leq p\leq z} \frac{(p-3)p^2}{(p-1)^3} 
\end{eqnarray}
we have the theorem. 
\end{proof}

Let $\#N$ be the number of prime factors in $N$. It is easy to see that $2^{\#N} < \sqrt{N}$ when $N$ is large. Thus, by taking $z = \sqrt{2N}$ in the previous theorem, we have $(N-2)\omega_P^1(N) > 2^{\#N-1}$ for large $N$. Also, we have $N \omega_P^1(N) > 4$ for large $N$. Numerical calculation shows that $N\geq \arabic{N0}$ will meet these conditions. 
\begin{theorem}
For $N\geq N_0\triangleq\arabic{N0}$, $(N-2)\omega_P^1(N) > 2^{\#N-1}$ and $N \omega_P^1(N) > 4$. 
\end{theorem}

We need this result in the next section.

\section{Even Number as a Sum of Two Distinct Primes}

We will prove in this section that, under {\bf UBH}, even number $2N$ is a sum of two distinct primes for $N \geq N_0 = \arabic{N0}$, more specifically, $S_P(N,N-2) \geq 1$ where $P$ is the product of all primes $\leq \sqrt{2N}$. First we introduce three lemmas. We assume $x>0$ in the following three lemmas. 
\begin{lemma}
Assume $p\mid P_{6N}$ and $d\mid P_{6pN}$. If $T_P^1(N,x) \leq T_{P_p}^1(N,x)$ and $T_P^1(N,x) \leq T_{P_d}^1(N,x)$, then $T_P^1(N,x) \leq T_{P_{dp}}^1(N,x) + t_{P_{dp}}^1(N)$.
\end{lemma} 
\begin{proof}
By assumption, $T_{P_d}^1(N,x) + T_{P_p}^1(N,x) - 2T_P^1(N,x)\geq 0$ and $t_P^1(N) = t_{P_p}^1(N) = t_{P_d}^1(N) = 0$. Thus, 
\begin{eqnarray}
S_{P_d}^1(N,x) + S_{P_p}^1(N,x) - 2S_P^1(N,x) \leq x\omega_{P_d}^1(N) + x\omega_{P_p}^1(N) - 2x\omega_P^1(N)
\end{eqnarray}
Since $S_{P_d}^1(N,x) + S_{P_p}^1(N,x) - S_P^1(N,x) = S_{P_{dp}}^1(N,x)$ by the inclusion-exclusion principle, then 
\begin{eqnarray}
S_{P_{dp}}^1(N,x) - S_P^1(N,x) \leq x\omega_{P_d}^1(N) + x\omega_{P_p}^1(N) - 2x\omega_P^1(N) 
\end{eqnarray} 
Now we have  
\begin{eqnarray}
\omega_{P_d}^1(N) + \omega_{P_p}^1(N) - \omega_P^1(N) & = & \omega_{P_{dp}}^1(N) \left(1-\frac{3}{p} + \prod_{q\mid d}\left(1-\frac{3}{q}\right) - \left(1-\frac{3}{p}\right)\prod_{q\mid d}\left(1-\frac{3}{q}\right)\right) \\
& = & \omega_{P_{dp}}^1(N) \left(1-\frac{3}{p} + \frac{3}{p}\prod_{q\mid d}\left(1-\frac{3}{q}\right)\right) \leq \omega_{P_{dp}}^1(N) 
\end{eqnarray}
Thus, $S_{P_{dp}}^1(N,x) - S_P^1(N,x) \leq x\omega_{P_{dp}}^1(N) - x\omega_P^1(N)$ and $T_P^1(N,x) \leq T_{P_{dp}}^1(N,x) + t_{P_{dp}}^1(N)$. 
\end{proof}

By this lemma, we get the following by induction method on $d\mid P_{6N}$: 
\begin{lemma} \label{Lemma_T_P_N_x}
If $T_P^1(N,x) \leq T_{P_p}^1(N,x)$ for each $p\mid P_{6N}$, then $T_P^1(N,x) \leq T_{P_d}^1(N,x)+ t_{P_{d}}^1(N)$ for any $d\mid P_{6N}$. 
\end{lemma}

\begin{lemma} \label{Lemma_sqrt_N}
Assume $c \mid 6$, $N$ is square-free and $c\perp N$. Let $P' = cN$. If $x$ is not integer, then $|T_{P'}^1(N,x)| \leq 2^{\#N-1}$ where $\#N$ is the number of prime factors in $N$. 
\end{lemma}
\begin{proof}
First $P'$ is squre-free and $N_{P'} = N$. For any $p\mid N$, we have $T_{P'}^1(N,x) = T_{P_p'}^1(N,x) - T_{P_p'}^1(\bar pN,x/p)$. Thus, we get the following after repeating use of this formula for all $p\mid N$: 
\begin{eqnarray}
T_{P'}^1(N,x) = \sum_{d\mid N} \mu(d) T_c^1\left(\bar dN, \frac{x}{d}\right)
\end{eqnarray}
Since $T_c^1(\bar dN, x/d) = T_1^1(\overline{cd}N, x/(cd))$ and 
\begin{eqnarray}
\left|T_1^1\left(\overline{cd}N, \frac{x}{cd}\right)\right| = \left|\left\{\frac{x}{cd}\right\}-\frac{1}{2}\right| \leq \frac{1}{2}
\end{eqnarray}
for any $d\mid N$, then 
\begin{eqnarray}
\left|T_{P'}^1(N,x)\right| = \left|\sum_{d\mid N} \mu(d) T_1^1\left(\overline{cd}N, \frac{x}{cd}\right) \right| \leq \sum_{d\mid N}\left|T_1^1\left(\overline{cd}N, \frac{x}{cd}\right)\right| \leq \sum_{d\mid N} \frac{1}{2} = 2^{\#N-1}
\end{eqnarray}
This completes the proof. 
\end{proof}

Let us start with the first upper bound hypothesis in terms of density. \\ \\
{\bf UBH:} First upper bound hypothesis on $S_P^p(N,x)$. Assume $N\geq N_0$ and $P$ is the product of all primes $\leq \sqrt{2N}$. For $p\mid P_{6N}$ and $N/2\leq x < N-1$, the following inequality holds: 
\begin{eqnarray} 
S_P^p(N,N\overline{P_p}P_p+x) - S_P^p(N,N\overline{P_p}P_p-x) \leq 3x\omega_P^p(N)
\end{eqnarray}
where $y = N\overline{P_p}P_p$ is constant. UBH can be given equivalently in terms of $T_P^p(N,x)$: \\ \\
{\bf UBH:} First upper bound hypothesis on $T_P^p(N,x)$. Assume $N\geq N_0$ and $P$ is the product of all primes $\leq \sqrt{2N}$. For $p\mid P_{6N}$ and $N/2\leq x < N-1$, the following inequality holds: 
\begin{eqnarray} \label{upper_bound_hypo_error} 
T_P^p(N,N\overline{P_p}P_p-x) - T_P^p(N,N\overline{P_p}P_p+x) \leq x\omega_P^p(N) 
\end{eqnarray}

\begin{theorem} 
Assume $N\geq N_0$ and $P$ is the product of all primes $\leq \sqrt{2N}$. If {\bf UBH} $(\ref{upper_bound_hypo_error})$ is true, then 
\begin{eqnarray} 
S_P(N,N-2) \geq 1
\end{eqnarray}
\end{theorem}
\begin{proof}
It is sufficient to prove $S_P^d(N,x) > 0$ for some $d\mid P_{6N}$ and some non-integer $x$ between $N/2$ and $N-1$. If there is $p\mid P_{6N}$ such that $T_P^p(N,x) < x\omega_P^p(N)$, then 
\begin{eqnarray}
S_P^p(N,x) = x\omega_P^p(N) - T_P^p(N,x) > 0
\end{eqnarray}
and the theorem is proved. Otherwise, we have $T_P^p(N,x) \geq x\omega_P^p(N)$ for each $p\mid P_{6N}$. From the first deduction formula for $T_P^1(N,x)$ and by {\bf UBH} $(\ref{upper_bound_hypo_error})$, we have 
\begin{eqnarray}
T_P^1(N,x) & = & T_{P_p}^1(N,x) - T_P^p(N,x) - T_P^p(N,N\overline{P_p}P_p+x) + T_P^p(N,N\overline{P_p}P_p-x) \\
& \leq & T_{P_p}^1(N,x) - T_P^p(N,x) + x\omega_P^p(N) \leq T_{P_p}^1(N,x) 
\end{eqnarray}

By Lemma \ref{Lemma_T_P_N_x}, we have $T_P^1(N,x) \leq T_{P_d}^1(N,x)+t_{P_d}^1(N)$ for any $d\mid P_{6N}$. Let $d' = P_{6N}$ and $P' = P_{d'}$, then $P' \mid 6N$ and $T_P^1(N,x) \leq T_{P'}^1(N,x)+t_{P'}^1(N)$ for any non-integer $x$ between $N/2$ and $N-1$. If $T_{P'}^1(N,x)+t_{P'}^1(N) \leq 2$, then for $N\geq N_0 = \arabic{N0}$,
\begin{eqnarray} 
S_P^1(N,x) = x\omega_P^1(N) - T_P^1(N,x) - t_{P'}^1(N) \geq x\omega_P^1(N) - 2 > 0
\end{eqnarray}
and the theorem is proved. Now we assume $P' \leq N$. Let $n = \lfloor N/P'\rfloor$ and $y = nP' + \frac{1}{2}$, then $N/2 \leq nP' \leq N$ and $y < N+1$. Since $T_{P'}^1(N,x)$ is periodic with period $P'$, then 
\begin{eqnarray}
T_{P'}^1(N,y) = T_{P'}^1\left(N,nP'+\frac{1}{2}\right) = T_{P'}^1\left(N,\frac{1}{2}\right) \leq \frac{1}{2}
\end{eqnarray}
and $T_{P'}^1(N,y)+t_{P'}^1(N)\leq 1$ since $t_{P'}^1(N) \leq \frac{1}{2}$. Now we have
\begin{eqnarray}
T_{P'}^1(N,y-2) & = & (y-2)\omega_{P'}^1(N) - t_{P'}^1(N) - S_{P'}^1(N,y-2) \\
& = & T_{P'}^1(N,y) - 2\omega_{P'}^1(N) + S_{P'}^1(N,y) - S_{P'}^1(N,y-2)
\end{eqnarray}
Since $d'$ is odd, $P'$ is even and $S_{P'}^1(N,y) - S_{P'}^1(N,y-2)\leq 1$, then $T_{P'}^1(N,y-2) \leq T_{P'}^1(N,y) + 1$ and $T_{P'}^1(N,y-2) + t_{P'}^1(N)\leq 2$. If $y < N-1$, then we take $x = y$; otherwise we take $x = y - 2 < N-1$. In either case we have $T_{P'}^1(N,x) +t_{P'}^1(N)\leq 2$ and the theorem is proved. Next we assume $P' > N\geq N_0$. In this case $P'$ has at least one odd prime factor $> 3$. Let $P'' = P'/2 = P_{2d'}$, then $P''\mid 3N$ and $P'' > N/2$. Since $W_{P'}^1(N)$ is symmetric in the middle between 1 and $P'$, then 
\begin{eqnarray}
S_{P'}^1(N,P'') = |W_{P'}^1(N)|\frac{P''}{P'} = \frac{1}{2}\prod_{p\mid P_3'}(p-1)
\end{eqnarray}

and $T_{P'}^1(N,P'') = 0$. If $P'' \leq N$, we take $y = P'' + \frac{1}{2}$, then $T_{P'}^1(N,y) + t_{P'}^1(N)\leq 1$. If $y < N-1$, then we take $x = y$; otherwise we take $x = y - 2$. In either case we have $T_{P'}^1(N,x) + t_{P'}^1(N) \leq 2$ and the theorem is proved. Finally, we assume $P'' > N$. In this case, $P'' = 3N$ or $P' = 3N$; or $P' = cN$ where $c$ is the index of $(P,N)$ and $c = 3$ or 6. Let $x = N - \frac{3}{2}$, then by Lemma \ref{Lemma_sqrt_N}, we have $T_{P'}^1(N,x) \leq 2^{\#N-1}$ and 
\begin{eqnarray}
S_P^1(N,x) = x\omega_P^1(N) - T_P^1(N,x) > (N-2)\omega_P^1(N) - 2^{\#N-1} > 0 
\end{eqnarray}
since $N\geq N_0 = \arabic{N0}$. That completes the proof. 
\end{proof}

We can further add one term on inequality $(\ref{upper_bound_hypo_error})$ as follows:
\begin{eqnarray} 
T_P^p(N,N\overline{P_p}P_p-x) - T_P^p(N,N\overline{P_p}P_p+x) \leq x\omega_P^p(N) + \theta(N)
\end{eqnarray}
where $\theta(N) = \sqrt{N}\omega_P^1(N)$. No example is found that this inequality fails for $N\geq 30,000,000$. 
\begin{lemma} 
Assume $x > 0$ and $P$ is the product of all primes $\leq \sqrt{2N}$. If $T_P^1(N,x) \leq T_{P_p}^1(N,x) + \theta(N)$ for each $p\mid P_{6N}$, then 
\begin{eqnarray}
T_P^1(N,x) \leq T_{P_d}^1(N,x) + \theta(N)\pi(\sqrt{2N}) + t_{P_d}^1(N)
\end{eqnarray}
for any $d\mid P_{6N}$. 
\end{lemma}

Under the condition of this lemma, we have, 
\begin{eqnarray}
S_P^1(N,N-2) & > & (N-2)\omega_P^1(N) - \theta(N)\pi(\sqrt{2N}) - 2^{\#N-1} \\
& \sim & \left(1 - \frac{2\sqrt{2}}{\log N}\right)N \omega_P^1(N)
\end{eqnarray}
and expect $S_P(N,N-2) > 0$ for $N$ larger than another fixed number.

\section{Twin Primes}

Goldbach conjecture says that for every $N\geq 4$, there is a pair of distinct primes $p$ and $q$ such that $q + p = 2N$. The generalized twin prime conjecture says that for every $N\geq 1$, there are infinitely many pairs of primes $p$ and $q$ such that $q - p = 2N$. This similarity gives the similar answer to both conjectures. Now we present the second upper bound hypothesis: \\ \\
{\bf UBH$'$:} Second upper bound hypothesis on $T_P^p(N,x)$. For given $N\geq 1$, there are infinitely many integers $M \geq 2N+1$ such that, for $M^2-7N \leq x < M^2-N$ and for each $p\mid P_{6N}$ where $P$ is the product of all primes $\leq M$, the following inequality holds: 
\begin{eqnarray} \label{upper_bound_hypo_twin} 
T_P^p(N,N\overline{P_p}P_p-x) - T_P^p(N,N\overline{P_p}P_p+x) \leq x\omega_P^p(N) 
\end{eqnarray}

\begin{lemma}
Let $M \geq 2N+1$ and $P$ the product of all primes $\leq M$. If $S_P(N,M^2-N) \geq 1$, then there exists $n\leq M^2-N$ such that both $n+N$ and $n-N$ are prime, and $n-N > M$. 
\end{lemma}
\begin{proof}
Since $S_P(N,M^2-N) \geq 1$, then there is $n$ between 1 and $M^2-N$ such that $(n+N)(n-N)\perp P$. Since $M\geq 2N+1$ and $n+N\perp P$, then $n+N\geq M+1\geq 2N+2$ and $n-N\geq 2$. Since $n+N \leq M^2$, then both $n+N$ and $n-N$ are prime. Since $n-N\perp P$, then $n-N > M$.
\end{proof}

\begin{lemma}
If there are infinitely many integers $M\geq 2N+1$ such that $S_P(N,M^2-N) \geq 1$ where $P$ is the product of all primes $\leq M$, then there are infinitely many pairs of primes $p$ and $q$ such that $q - p = 2N$. 
\end{lemma}
\begin{proof}
We choose $M_1\geq 2N+1$ and $M_{m+1} \geq M_m^2$ for each $m\geq 1$. By the assumption, we have $S_P(N,M_m^2-N) \geq 1$ where $P$ is the product of all primes $\leq M_m$. Thus, there is $n_m$ between 1 and $M_m^2-N$ such that, by the previous lemma, both $p_m = n_m-N$ and $q_m = n_m+N$ are prime and $q_m-p_m = 2N$. Since 
\begin{eqnarray}
q_m = n_m+N \leq M_m^2 \leq M_{m+1} < n_{m+1} - N = p_{m+1} < q_{m+1}
\end{eqnarray}
then $q_m < q_{m+1}$. Thus, there are infinitely many pairs of primes $p$ and $q$ such that $q - p = 2N$.
\end{proof}

\begin{theorem} 
If {\bf UBH$'$} $(\ref{upper_bound_hypo_twin})$ is true for given $N\geq 1$, then there are infinitely many pairs of primes $p$ and $q$ such that $q - p = 2N$. 
\end{theorem}
\begin{proof}
Let $M$ be the one of integers in {\bf UBH$'$} and $P$ the product of all primes $\leq M$. It is sufficient to prove $S_P(N,M^2-N) \geq 1$ due to the previous lemma. Assume $x$ is not integer and $M^2-7N < x < M^2-N$. If there is $p\mid P_{6N}$ such that $T_P^p(N,x) < x\omega_P^p(N)$, then 
\begin{eqnarray}
S_P^p(N,x) = x\omega_P^p(N) - T_P^p(N,x) > 0
\end{eqnarray}
and the theorem is proved. Otherwise, for each $p\mid P_{6N}$, we have $T_P^p(N,x) \geq x\omega_P^p(N)$. From the first deduction formula for $T_P^1(N,x)$ and by {\bf UBH$'$} $(\ref{upper_bound_hypo_twin})$, we have
\begin{eqnarray}
T_P^1(N,x) & = & T_{P_p}^1(N,x) - T_P^p(N,x) - T_P^p(N,N\overline{P_p}P_p+x) + T_P^p(N,N\overline{P_p}P_p-x) \\
& \leq & T_{P_p}^1(N,x) - T_P^p(N,x) + x\omega_P^p(N) \leq T_{P_p}^1(N,x)
\end{eqnarray}

By Lemma \ref{Lemma_T_P_N_x}, we have $T_P^1(N,x) \leq T_{P_d}^1(N,x)+t_{P_d}^1(N)$ for any $d\mid P_{6N}$. Let $d' = P_{6N}$ and $P' = P_{d'}$, then $P' \mid 6N$. Thus, there is $n$ such that $M^2-7N \leq nP' < M^2-N$. Let $x = nP'+\frac{1}{2}$, then $T_{P'}^1(N,x)+t_{P_d}^1(N) \leq 1$ and 
\begin{eqnarray}
S_P^1(N,x) = x\omega_P^1(N) - T_P^1(N,x) \geq (M^2-7N)\omega_P^1(N) - 1 > 0
\end{eqnarray}
That completes the proof. 
\end{proof}

We can further add one term on inequality $(\ref{upper_bound_hypo_twin})$ as follows:
\begin{eqnarray} \label{upper_bound_hypo_error_more} 
T_P^p(N,N\overline{P_p}P_p-x) - T_P^p(N,N\overline{P_p}P_p+x) \leq x\omega_P^p(N) + \theta'(M)
\end{eqnarray}
where $\theta'(M) = M \omega_P^1(N)$, and have the following lemma: 
\begin{lemma} 
Assume $x>0$ and $P$ is the product of primes $\leq M$. If $T_P^1(N,x) \leq T_{P_p}^1(N,x) + \theta'(M)$ for each $p\mid P_{6N}$, then 
\begin{eqnarray}
T_P^1(N,x) \leq T_{P_d}^1(N,x) + \theta'(M)\pi(M)+t_{P_d}^1(N)
\end{eqnarray}
for any $d\mid P_{6N}$. 
\end{lemma}

Under the condition of this lemma, we have 
\begin{eqnarray}
S_P^1(N,M^2-N) & > & (M^2-7N)\omega_P^1(N) - \theta'(M)\pi(M) - 1 \\
& \sim & \left(1 - \frac{1}{\log M}\right) M^2\omega_P^1(N)
\end{eqnarray}
and expect $S_P^1(N,M^2-N) > 0$ for $M$ larger than a fixed number. 

\section*{Acknowledgment}

I would like to thank Dr David Platt of University of Bristol for his pointing out that {\bf UBH} fails when $N = 400$ and $p = 23$. Author believes {\bf UBH} is valid for $N$ large enough. Numerical calculation shows {\bf UBH} $(\ref{upper_bound_hypo_first})$ is valid for $N$ between $100,000,000$ and $102,000,000$.

\end{document}